\documentclass[reqno,12pt]{amsart}

\usepackage[colorlinks=true]{hyperref}
\hypersetup{urlcolor=blue, citecolor=red}

\usepackage{amssymb,graphicx,color}
\usepackage[latin1]{inputenc}

\theoremstyle{plain}
\newtheorem{thm}{Theorem}[section]

\theoremstyle{definition}
\newtheorem{rmk}{Remark}[section]

\newcommand{\NN}{{\mathbb{N}}}
\newcommand{\ZZ}{{\mathbb{Z}}}
\newcommand{\RR}{{\mathbb{R}}}

\newcommand{\bfu}{\mathbf{u}}
\newcommand{\bfv}{\mathbf{v}}
\newcommand{\bfw}{\mathbf{w}}
\newcommand{\bfe}{\mathbf{e}}

\newcommand{\bfg}{\mathbf{g}}

\newcommand{\bfx}{\mathbf{x}}

\newcommand{\bfU}{\mathbf{U}}
\newcommand{\bfF}{\mathbf{F}}

\newcommand{\bfzero}{\mathbf{0}}

\newcommand{\bfnabla}{\boldsymbol{\nabla}}

\newcommand{\Ccal}{{\mathcal C}}

\newcommand{\Ical}{{\mathcal I}}

\newcommand{\Kcal}{{\mathcal K}}

\newcommand{\Vcal}{{\mathcal V}}

\newcommand{\Pcal}{{\mathcal P}}

\newcommand{\Rey}{{\text{\rm Re}}}

\newcommand{\rmd}{{\text{\rm d}}}

\newcommand{\average}[1]{\langle{#1}\rangle}

\newcommand{\dinner}[1]{(\!({#1})\!)}

\newcommand{\dontshow}[1]{}

\begin{document}
\numberwithin{equation}{section}

\makeatletter
\@namedef{subjclassname@2020}{\textup{2020} Mathematics Subject Classification}
\makeatother


\title[Energy cascade in free-shear 3D flows]{Energy cascade for cross-shear length scales in free-shear three-dimensional incompressible viscous flows}

\author[R. Rosa]{Ricardo M. S. Rosa}

\address[R. Rosa]{Instituto de Matem\'atica, Universidade Federal do Rio de Janeiro, Brazil}

\email[R. Rosa]{rrosa@im.ufrj.br}

\date{\today}

\dedicatory{To my mentor, collaborator and friend, Roger Temam, on the occasion of his 85th birthday.}

\thanks{This work was partly supported by Coordena\c{c}\~ao de Aperfei\c{c}oamento de Pessoal de N\'{\i}vel Superior (CAPES), Brasil, grant 001.}



\begin{abstract}
    The phenomenon of energy cascade is addressed in the case of free-shear flows, modeled with the equations for incompressible Newtonian fluids with mixed periodic and free-slip boundary conditions driven by an imposed mean shear profile. The rigorous results are proved for ensemble averages with respect to stationary statistical solutions in the sense of Foias and Prodi. The flow is decomposed into high and low horizontal wavenumber components, corresponding to scales perpendicular to the mean shear gradient, and the corresponding energy budget balances viscous dissipation, production by the mean shear, and transfer between scales. Exploiting the orthogonality of horizontal Fourier modes to control the shear-production term, we look for a cascade in which energy flows from the large scales, where it is fed by the mean shear, down to the small scales, where viscosity dissipates it. For any given wavenumber, we define an associated horizontal Taylor wavenumber of the high-wavenumber band above the given wavenumber. We prove two energy cascade results, for the energy flux and for a restricted energy flux that discounts possible energy loss to flow singularities, both valid for wavenumbers whose associated horizontal Taylor wavenumber lies above the viscous shear wavenumber and whose low-wavenumber energy dissipation rate is negligible compared to the total energy dissipation rate. On heuristic grounds, these rigorous conditions correspond precisely to an energy cascade in the range from the classical Corrsin scale down to the Kolmogorov dissipation scale, as expected in the conventional theory of turbulent shear flows.
\end{abstract}

\maketitle

\tableofcontents

\section{Introduction}

We consider statistically stationary free-shear flows and look for conditions for the existence of an inertial range in which energy cascades from larger to smaller cross-shear scales, i.e. length scales in the two-dimensional plane perpendicular to the shear direction.

The shear flow model we consider is the following incompressible Navier--Stokes system, driven by a given shear flow $\bfU = (U(z),0,0)$, and complemented with mixed periodic and free-slip boundary conditions and with the streamwise and spanwise velocity components having zero horizontal average over the cross-shear planes:
\begin{equation}
  \label{nsefluctuation}
  \begin{cases}
    \displaystyle \frac{\partial \bfu}{\partial t} - \nu \Delta \bfu + (\bfu\cdot\bfnabla)\bfu + \bfnabla p = - U(z)\frac{\partial\bfu}{\partial x} - w U'(z) \bfe_x,  \\
    \bfnabla \cdot \bfu = 0, \\
    \bfu(\bfx + L_x \bfe_x) = \bfu(\bfx), \quad p(\bfx + L_x \bfe_x) = p(\bfx), \\
    \bfu(\bfx + L_y \bfe_y) = \bfu(\bfx), \quad p(\bfx + L_y \bfe_y) = p(\bfx), \\
    \displaystyle \frac{\partial u}{\partial z}(x,y,\pm h/2) = 0, \;\frac{\partial v}{\partial z}(x,y,\pm h/2) = 0, \;w(x,y,\pm h/2) = 0, \\
    \displaystyle \int_0^{L_x}\int_0^{L_y} u(x,y,z)\;\mathrm{d}x\,\mathrm{d}y = 0, \quad \int_0^{L_x}\int_0^{L_y} v(x,y,z)\;\mathrm{d}x\,\mathrm{d}y = 0.
  \end{cases}
\end{equation}
The unknowns are the velocity field $\bfu = (u, v, w)$ and the kinematic pressure $p,$ which are functions of the time variable $t$ and the space variable $\bfx = (x, y, z).$ The canonical basis in the $xyz$ space is denoted by $\{\bfe_x, \bfe_y, \bfe_z\}$. The parameter $\nu > 0$ is the kinematic viscosity of the fluid. The spatial domain is characterized by the periods $L_x$, $L_y > 0$ in the $x$ and $y$ directions, respectively, and by the characteristic height $h > 0$ of the shear flow. The imposed shear profile $U=U(z)$ is assumed to have uniformly bounded derivative.

The interpretation of the system \eqref{nsefluctuation} is that we are looking at the fluctuation $\bfu$ from a given shear flow $\bfU$, which is taken to mimic a turbulent mean shear flow. The system is driven by the shear mechanism, which is responsible for the generation of small scales in the flow, and which is represented by the terms on the right-hand side of the momentum equation. The term $-U(z)\partial_x\bfu$ is an advection term, and $-w U'(z) \bfe_x$ is a shear-production term. The presence of these terms is one aspect that makes the problem different from the previous works on energy cascades in the Navier--Stokes equations, which were all driven purely by body forces.

The model \eqref{nsefluctuation} can be formally derived by writing the Navier--Stokes equations for the total velocity field $\bfU + \bfu$, under the assumption that the mean flow is a stationary and horizontally homogeneous parallel flow, together with a suitable ergodic hypothesis and closure assumptions on the Reynolds stresses (see Section \ref{secmodel} for the details).

The mixed periodic and free-slip boundary conditions are imposed on the fluctuation field $\bfu$ directly. The mean shear profile $\bfU$ is prescribed as a driving term, encoding the conditions that sustain the shear without modeling them explicitly (e.g. the effect of inflow and outflow conditions or non-homogeneous slip boundary conditions). The pressure $p,$ on the other hand, is the combination of the mean and fluctuation kinematic pressures.

The existence of energy cascades in turbulent flows is well-documented in experiments and is backed by phenomenological arguments \cite{batchelor1953, Champagne1970, tennekeslumley1972, brownroshko1974, hinze1975, michalke1964, moninyaglom1975, Tavoularis1981, RogersMoinReynolds1986, frisch1995, lesieur1997, pope2000, jimenez2012, zadpk2022}. Our aim is to obtain conditions for the cascade from \emph{first principles}, based on the Navier--Stokes equations. More precisely, we look for conditions pertaining to stationary statistical solutions in the sense of \cite{foias72,foias73,vishikfursikov88,fmrt2001a,frt2010b}, which includes, in particular, the (generalized) asymptotic limit of time-averages of Leray-Hopf weak solutions of the Navier--Stokes system (see \cite{foiastemam1975, foiasprodi76, fmrt2001a, frt2010b, frt2015, frt2019}). Such a statistical solution is directly connected to the notion of ensemble average. So, this is a rigorous framework to obtain conditions for the cascade in average.

In our setting, we consider ensemble averages $\average{\cdot}$ formally defined as averages with respect to a stationary statistical solution in the sense of Foias and Prodi (Section \ref{secsss}). We split the velocity field $\bfu$ into components with modes with horizontal wavenumber smaller and higher than a threshold horizontal wavenumber $\bar\kappa$, i.e. $\bfu = \bfu_{0, \bar\kappa} + \bfu_{\bar\kappa, \infty}$ (see Section \ref{secspecdecomposition}). Then we obtain the mean energy-budget inequalities (see \eqref{meanenergybudgetineq} and \eqref{meanenergybudgetineqfluct} in Section \ref{secenergybudget})
\[
  \nu\kappa_1^3\average{\|\bfnabla\bfu_{\bar\kappa,\infty}\|_{L^2}^2} \leq - \kappa_1^3\average{\dinner{(\bfu\cdot\bfnabla)\bfu_{0, \bar\kappa}, \bfu_{\bar\kappa, \infty}}_{L^2}}- \kappa_1^3\average{\dinner{w U'(z) \bfe_x,\bfu_{\bar\kappa,\infty}}_{L^2}},
\]
where $\kappa_1$ is the smallest possible wavenumber. The first term accounts for the mean energy dissipation rate within the higher horizontal modes, while the second term represents the flux of energy towards the higher modes, and the last term represents the production of energy by the mean shear flow within those modes. The flux term is denoted by $\overset{\rightarrow}{\Pi}_{\bar\kappa}$ and is the object we are most interested in. Our aim is to obtain conditions for this flux to be nearly constant and approximately equal to the (total) mean energy dissipation rate $\epsilon = \nu\kappa_1^3\average{\|\bfnabla\bfu\|_{L^2}^2}$. 

By estimating the production term, we prove in Theorem \ref{thmenergyfluxlowerbound} (see also Remark \ref{rmklowerboundcascade}) that 
\[
    \overset{\rightarrow}{\Pi}_{\bar\kappa} \gtrapprox \epsilon,
\]
for wavenumbers $\bar\kappa$ such that
\[
    \frac{{\kappa_S}^2}{2} \ll \Kcal_{\bar\kappa,\infty}^2 \quad \text{and} \quad \epsilon_{0, \bar\kappa} \ll \epsilon,
\]
where ${\kappa_S} = (S/\nu)^{1/2}$ is the viscous shear wavenumber, $S = \|U'\|_{L^\infty}$ measures the shear {rate,} $\epsilon_{0, \bar\kappa} = \nu\kappa_1^3\average{\|\bfnabla\bfu_{0, \bar\kappa}\|_{L^2}^2}$ is the energy dissipation rate within the lower horizontal wavenumbers, and 
\[
    \Kcal_{\bar\kappa,\infty} = \left(\frac{\average{\|\bfnabla_{x,y}\bfu_{\bar\kappa,\infty}\|_{L^2}^2}}{\average{\|\bfu_{\bar\kappa,\infty}\|_{L^2}^2}}\right)^{1/2}
\]
is a wavenumber depending on the threshold wavenumber $\bar\kappa.$

In Theorem \ref{thmrestrictedenergyfluxbounds}, we consider a restricted energy flux $\overset{\rightarrow}{\Pi}_{\bar\kappa}^*$, that discounts the possible energy loss due to a loss of regularity of the stationary statistical solutions, i.e. discounting a ``flux to modes at infinity'', and for which we have an equality in the energy budget (see equation \eqref{meanenergybudgeteqrestrited}), and prove that, for the same range of wavenumber $\bar\kappa$ as above, we have the energy cascade
\[
    \overset{\rightarrow}{\Pi}_{\bar\kappa}^* \approx \epsilon.
\]

As discussed in the concluding remarks, more specifically in Sections \ref{concludingrmkinertialrange} and \ref{secconcludingcomparison}, the conditions on the inertial range for the cascade are not explicit in the wavenumber $\bar\kappa,$ unlike in previous works on rigorous results for turbulent cascades. The explicit condition obtained by bounding the quantity $\Kcal_{\bar\kappa,\infty}$ from below by $\bar\kappa$ is successfully used in previous works, but here, due to the shear production of energy, the range obtained with the lower bound is not realistic since it includes ranges where the shear forces are significant, and the flow on those scales is not isotropic. In shear flows, for sufficiently large shear rates, the inertial range is supposed to occur at wavenumbers higher than the Taylor wavenumber used as an upper bound for the cascade in previous works. More precisely, it is expected instead that the inertial range lies between the Corrsin scale, defined by the wavenumber $\kappa_C = S^{3/2} / \epsilon^{1/2},$ and the Kolmogorov dissipation scale, defined by the wavenumber $\kappa_\eta = (\epsilon/\nu^3)^{1/4}$. By retaining $\Kcal_{\bar\kappa,\infty}$ instead of using its lower bound $\bar\kappa$ we improve considerably the rigorous estimate for the inertial range, at the cost of not having an explicit range. It turns out that, assuming the $-5/3$ Kolmogorov spectrum within the Corrsin and the Kolmogorov scale, the rigorous conditions we obtained read $\kappa_S^3/\kappa_\eta^2 \ll \bar\kappa \ll \kappa_\eta$, with the lower bound being precisely $\kappa_S^3/\kappa_\eta^2 = \kappa_C$, yielding exactly the expected inertial range from the conventional theory of shear turbulence. Hence, in contrast to the standard rigorous treatments of body-forced turbulence mentioned below, where the lower edge of the inertial range is essentially dictated by the forcing scale, a distinctive feature of our work is that our implicit lower bound is genuinely driven by the shear.

Similar results with conditions for the energy cascade were first proved in the case of fully-periodic three-dimensional flows with body forces \cite{fmrt2001c,fjmrt05}, followed by results on the enstrophy cascade in two-dimensional flows \cite{fjmr02}. Subsequent works considered the energy cascade in the physical space \cite{dascaliucgrujic2011, dascaliucgrujic2012a, dascaliucgrujic2012b}, as well as cascades in spectral space for passive scalars \cite{jollywirosoetisno2018}, and vertically-averaged convection \cite{balciisenbergjolly2018}.

The periodic channel driven by a pressure gradient was considered in \cite{rrt08} (after \cite{constantindoering95}; see also \cite{howard72}), but with the different aim of improving the estimate for the mean energy dissipation rate. The energy cascade was investigated, but not obtained. The difficulty in adapting the result from \cite{fmrt2001c,fjmrt05} is that kinetic energy is generated at arbitrarily small scales, with a slow decay rate with respect to the scales. This is due to the presence of a boundary layer near the walls, where the shear is large, and which is responsible for the generation of small scales.

In the current article, we analyse instead the above mathematical model of free-shear flow, where we look at the fluctuations from a given shear flow imposed to mimic a turbulent mean flow. This drives the system in a way different from the previous problems proving energy cascade results, which were driven by body forces. Two of those, namely \cite{balciisenbergjolly2018} and \cite{jollywirosoetisno2018}, also add a convection-driven mechanism in the passive scalar term considered along the fluid flow, but the velocity fields are driven purely by body forces, and the argument to prove the cascade for the passive scalar is of a different nature. Our results look specifically at the energy flux of the fluctuation term of the velocity field driven by a direct shear mechanism. So, this is indeed the first work to obtain rigorous energy cascade results with a force that is not purely a body force.

The other significant difference is that we look at the transfer of energy between the horizontal scales of motion, splitting the fluctuation velocity field between large and small horizontal scales, i.e. looking at the horizontal wavenumbers. Due to the shear gradient along the vertical direction, the splitting in the horizontal direction not only facilitates the analysis but captures the cascade at scales that might not be fully isotropic. The splitting with respect to the full wavenumber is also possible and complements the current work. For that, the conditions for the energy cascade are more stringent, requiring four conditions instead of two. The results for the full wavenumber splitting will be presented elsewhere.

The paper is organized as follows: Section \ref{secframework} sets up the framework we work on and gives the main definitions; Section \ref{secenergybudgets} proves the main energy-budget relations and the main estimates of the separate terms in the relations; Section \ref{secenergycascade} proves the main results about the energy cascade; and Section \ref{secconcludingrmks} makes some concluding remarks.

\section{Framework}
\label{secframework}

In this section we give more details on the model, set up the functional analysis framework, recall the definition of stationary statistical solutions for the ensemble averages, and introduce the physical quantities to be studied.

\subsection{Model motivation}
\label{secmodel}

In a given flow with a pressure field $q$, decomposing the velocity field as a sum $\bfv = \bfU + \tilde\bfu$ of a mean velocity field $\bfU = \average{\bfv}$, and a fluctuation component $\tilde\bfu = \bfv - \bfU$, and assuming \emph{the flow is statistically stationary} yield the momentum equation
\begin{equation}
    \label{nsemomentumdecomposition}
    \frac{\partial \tilde\bfu}{\partial t} + (\tilde\bfu \cdot \bfnabla) \tilde\bfu + (\bfU \cdot \bfnabla) \tilde\bfu + (\tilde\bfu \cdot \bfnabla) \bfU + (\bfU \cdot \bfnabla) \bfU + \bfnabla q = \nu \Delta \tilde\bfu + \nu \Delta \bfU.
\end{equation}
At this point, the mean flow is unknown explicitly, but it satisfies the stationary Reynolds-averaged Navier--Stokes equations, which are obtained by taking the ensemble average $\average{\cdot}$ of the above equation. This yields
\begin{equation}
    \label{nsereynoldsaverage}
    (\bfU \cdot \bfnabla) \bfU + \bfnabla P = \nu \Delta \bfU - \bfnabla \cdot \average{\tilde\bfu \otimes \tilde\bfu},
\end{equation}
where $P = \average{q}$ is the mean pressure, and $\average{\tilde\bfu \otimes \tilde\bfu}$ is the Reynolds stress tensor, and we denote the pressure fluctuation by $\tilde p= q - P.$ The mean flow $\bfU$ is then determined by the above equation, together with the incompressibility condition $\bfnabla \cdot \bfU = 0$ and with the mean of whatever the boundary conditions are for the original flow.

Assuming further that \emph{the mean flow is horizontally homogeneous and parallel to the $x$ direction,} we deduce that its velocity field must be of the form $\bfU = (U(z), 0, 0) = U(z) \bfe_x$. In this case, the term $(\bfU \cdot \bfnabla) \bfU$ drops out of the equations. Assuming also that \emph{the mean pressure has no horizontal gradient}, we deduce that $P=P(z)$ and that the Reynolds stress tensor $\average{\tilde\bfu\otimes\tilde\bfu}(z)$ is also a function of $z$ alone. Then, we obtain the momentum equation
\begin{equation}
    \label{nsemomentumdecomposition2}
    \frac{\partial \tilde\bfu}{\partial t} + (\tilde\bfu \cdot \bfnabla) \tilde\bfu + (\bfU \cdot \bfnabla) \tilde\bfu + (\tilde\bfu \cdot \bfnabla) \bfU + \bfnabla \tilde p + P'(z)\bfe_z = \nu \Delta \tilde\bfu + \nu U''(z)\bfe_x
\end{equation}
and the Reynolds equation
\begin{equation}
    \label{nsereynoldsaverage2}
    P'(z)\bfe_z = \nu U''(z)\bfe_x - \frac{\partial}{\partial z} \average{\tilde w \tilde\bfu}.
\end{equation}
Equations \eqref{nsemomentumdecomposition2} and \eqref{nsereynoldsaverage2} form a system of partially coupled equations which is not closed since the Reynolds stress tensor is unknown. An approach to modeling this type of shear flow is to \emph{drop the Reynolds equation \eqref{nsereynoldsaverage2}} and instead \emph{prescribe the mean flow} in equation \eqref{nsemomentumdecomposition2}; the prescribed mean flow then drives the fluctuation field $\tilde\bfu$. With the mean flow prescribed, \eqref{nsereynoldsaverage2} can be seen as a \emph{closure assumption} on the ensemble-averaged Reynolds stresses.

Assuming \emph{the prescribed velocity field $\bfU$ is of the form $\bfU = (U(z), 0, 0) = U(z) \bfe_x$ and the prescribed pressure $P$ is of the form $P=P(z)$} leads to the momentum equation
\begin{equation}
    \label{nsefluctuationprime}
    \frac{\partial \tilde\bfu}{\partial t}  + (\tilde\bfu\cdot\bfnabla)\tilde\bfu + U(z)\frac{\partial \tilde\bfu}{\partial x} + \tilde w U'(z) \bfe_x + \bfnabla q = \nu \Delta \tilde\bfu + \nu U''(z)\bfe_x,
\end{equation}
where $q = P + \tilde p$.

In the case that the shear rate $U'(z)$ imposed by the mean profile $U=U(z)$ is localized in a certain region $|z| \lesssim \delta$, where $\delta$ is a shear-core thickness, with $U'(z)$ nearly vanishing away from this region, it is expected that the fluctuations become negligible far away, where the mean flow carries no tangential viscous stress. This suggests, as another modeling assumption, \emph{confining the flow to a region $|z| < h/2$}, for a given $h > 0$, and \emph{imposing on $\tilde\bfu$ a free-slip condition at the boundaries $z=\pm h/2$}, which does not introduce any tangential viscous stress at the boundary. We also assume that the fluctuation $\tilde\bfu$ is \emph{periodic in the streamwise and spanwise directions, with respective periods $L_x$ and $L_y$.}

The final modeling assumption concerns the horizontal average of the fluctuation. Taking the spatial average of \eqref{nsefluctuationprime} over the cross-shear planes, one finds that the horizontal average of $\tilde\bfu$ is driven by the body-type force $\nu U''(z)\bfe_x$ and by the divergence of the instantaneous Reynolds stress, so it is not expected to vanish in general. Its ensemble average, however, vanishes, since $\langle\tilde\bfu\rangle = 0$. Under the ergodic hypothesis, averages over the homogeneous cross-shear directions approximate the ensemble averages, so the horizontal averages are expected to remain small at each instant of time when the periods $L_x$ and $L_y$ are large compared with the integral scale of the flow. Based on that, we assume that the horizontal average of the velocity field $\bfv$ is held fixed at the prescribed profile $U(z)$ at all times, i.e. that \emph{the horizontal average of the fluctuation $\tilde\bfu$ vanishes identically}. Under this assumption, \eqref{nsefluctuationprime} splits into two parts. Its component with zero horizontal average, in which the body-type force $\nu U''(z)\bfe_x$ does not appear, gives the momentum equation in \eqref{nsefluctuation} for the fluctuation, which we now denote simply by $\bfu$. Its horizontal average reduces to the Reynolds-like equations
\[
  \frac{\partial}{\partial z}\overline{\tilde w\tilde u} = \nu U''(z), \qquad \frac{\partial}{\partial z}\overline{\tilde w\tilde v} = 0,
\]
where the overline denotes the average over the cross-shear planes. With the profile $U=U(z)$ given a~priori, the above equations amount to a \emph{closure assumption} on the horizontally averaged Reynolds stresses and are dropped from the final system, similarly to dropping the Reynolds equation \eqref{nsereynoldsaverage2} when prescribing the mean flow. Thus, we are left with the equations \eqref{nsefluctuation} as a closed system for the fluctuation $\bfu$, subject to the zero-horizontal-average and free-slip conditions.

Classical free-shear flows suggest natural forms for the profile $U(z)$ (see e.g. \cite{michalke1964,Champagne1970,brownroshko1974,Tavoularis1981,RogersMoinReynolds1986,pope2000,lesieur1997}). For a \textit{mixing layer}, the standard choice is the hyperbolic tangent profile
\[
    U(z) = \frac{U_1 + U_2}{2} + \frac{U_1 - U_2}{2}\tanh\!\left(\frac{2z}{\delta}\right),
\]
where $U_1, U_2$ are the two free-stream velocities and $\delta$ is a characteristic shear-core thickness. This profile is well-supported experimentally and has a single inflection point at $z=0$, consistent with Rayleigh's instability criterion. For a \textit{planar jet}, one uses either
\[
    U(z) = U_0\,\operatorname{sech}^2\!\left(\frac{z}{\delta}\right) \quad \text{or} \quad U(z) = U_0\exp\!\left(-\frac{z^2}{2\delta^2}\right),
\]
and for a \textit{wake},
\[
    U(z) = U_\infty - U_d\exp\!\left(-\frac{z^2}{2\delta^2}\right),
\]
where $U_\infty$ is the free-stream velocity and $U_d > 0$ is the velocity deficit.

These classical mean profiles are considered in the whole space, with $U'(z)$ decaying at infinity, as in the motivation leading to our free-slip model. In our results, we do not commit to any particular profile, and all results are stated in terms of $U(z)$ and its derivatives. We only assume it is a smooth function of $z$ whose derivative $U'(z)$ is uniformly bounded. In particular, $U'$ need not vanish at $z = \pm h/2$; the decay of $U'$ only serves to motivate the model.

Another important profile, whose shear rate does not decay at infinity, is the linear profile $U(z) = Sz$, with constant shear rate $S$. In \cite{Champagne1970,Tavoularis1981}, experimental setups were carefully designed to approximate this profile in laboratory flows, while DNS simulations were performed with this profile imposed on the system in \cite{RogersMoinReynolds1986,SekimotoDongJimenez2016}; see Section \ref{concludingrmkinertialrange} for more details. In this case, the more natural boundary conditions are shear-periodic conditions in the $z$-direction, as considered in \cite{RogersMoinReynolds1986,SekimotoDongJimenez2016} (and implemented with Rogallo's periodic remeshing method \cite{Rogallo1981}), which better models homogeneity in all directions. On the other hand, our free-slip condition on the fluctuation field is a modeling compromise which blocks wall-normal transport in a thin layer near $z=\pm h/2$, and thereby locally suppresses shear production. Since the cascade result is a global statement, and the inertial-range dynamics are an interior small-scale phenomenon, with the boundary effect being confined to a region of relatively small volume, the free-slip and shear-periodic models are likely to share similar interior cascades. The proper mathematical treatment of the shear-periodic boundary conditions needs substantial work and is not pursued here.

\subsection{The flow domain}

We consider the flow domain
\[
    \Omega=\left(0,L_x\right)\times\left(0,L_y\right)\times\left(-h/2,h/2\right).
\]
In this geometry, the $x, y,$ and $z$ directions are called, respectively, the streamwise, the spanwise, and the shear-normal directions. The $x$ and $y$ directions are called cross-shear directions, and planes parallel to the $xy$ plane are called cross-shear planes. The boundaries at $z = \pm h/2$ are called the free-slip boundaries, and the boundaries at $x = 0, L_x$ and $y = 0, L_y$ are called the periodic boundaries. The free-slip boundary conditions prescribe no normal flow and vanishing tangential stress at $z = \pm h/2$, as given explicitly in \eqref{nsefluctuation}.

\subsection{Functional spaces}

As usual in the theory of weak solutions of the Navier--Stokes equations \cite{leray1934, lady63, temam1984, temam1995, constantinfoias89}, we consider the phase spaces to be closures of the space of smooth vector fields satisfying the essential boundary conditions. In our case, we consider the space $\Vcal$ of restrictions to $\Omega$ of divergence-free vector fields $\bfv=(u, v, w)$ in $\Ccal^\infty(\RR^2 \times (-h/2, h/2))$ which are $L_x$- and $L_y$- periodic in the $x$ and $y$ directions and such that $\partial_z u,$ $\partial_z v$ and $w$ have support bounded away from the boundaries $z = \pm h/2.$ The spaces $H$ and $V$ are the closures of $\Vcal$ with respect to the $L^2(\Omega)^3$ and $H^1(\Omega)^3$ norms, respectively.

With the assumption of zero horizontal average of the streamwise and spanwise components of the velocity field in mind, we also define the spaces
\[ 
    \dot H = \left\{\bfu \in H; \; \int_0^{L_x}\int_0^{L_y} (u(x,y,z), v(x,y,z))\;\mathrm{d}x\,\mathrm{d}y = (0, 0), \text{ for a.e.}\; z\right\},
\]
and
\[
    \dot V = V \cap \dot H.
\]

The inner product and norm in $H$ and in $\dot H$ are those inherited from $L^2$:
\[ \dinner{\bfu,\bfv}_{L^2} = \int_\Omega \bfu \cdot \bfv \;\rmd \bfx, \quad \|\bfu\|_{L^2} = \left(\int_\Omega |\bfu|^2 \;\rmd \bfx\right)^{1/2},
\]
with $|\bfu|^2 = u^2 + v^2 + w^2$. The inner product and norm in $V$ are those inherited from $H^1$, but for $\dot V,$ we consider
\[ \dinner{\bfnabla\bfu,\bfnabla\bfv}_{L^2} = \int_\Omega \bfnabla\bfu : \bfnabla \bfv \;\rmd\bfx, \quad \|\bfnabla\bfu\|_{L^2} = \left(\int_\Omega |\bfnabla\bfu|^2\;\rmd\bfx\right)^{1/2},
\]
where $\bfnabla\bfu$ is the velocity gradient tensor. Since the domain is bounded in the shear-normal direction and constant vector fields are excluded from $\dot V$ by the condition of zero horizontal average of the streamwise and spanwise components and by the condition that the shear-normal component vanishes on the free-slip boundaries, the Poincar\'e inequality holds in $\dot V$, and the above norm and inner product are equivalent to the usual $H^1$-norm and inner product.

Changing the notation momentarily and writing $\bfu=(u_i)_i$ and $\bfx=(x_j)_j$, we have $\bfnabla\bfu = (\partial_{x_j}u_i)_{i,j}$, $\bfnabla\bfu : \bfnabla\bfv = \sum_{i,j} \partial_{x_j}u_i\partial_{x_j}v_i$, and $|\bfnabla\bfu|^2 = \sum_{i,j} |\partial_{x_j}u_i|^2$. We also use the $L^\infty(\Omega)$ norm, denoted simply by $\|\bfu\|_{L^\infty},$ for vector fields $\bfu=\bfu(\bfx),$ and $\|U\|_{L^\infty},$ for scalar fields $U=U(z)$.

The Leray-Helmholtz projector $P_{\textrm{LH}}$ is the orthogonal projection from $L^2(\Omega)^3$ onto $H$, and we have the usual decomposition $L^2(\Omega)^3 = H \oplus H^\perp,$ where $H^\perp$ is the orthogonal complement of $H$ in $L^2(\Omega)^3$. 

The orthogonal projector onto the space of velocity fields with horizontally averaged streamwise and spanwise components and zero shear-normal component is given by
\[
    Q_h \bfu = \frac{1}{L_x L_y}\int_0^{L_x}\int_0^{L_y} (u(x,y,z), v(x,y,z), 0)\;\mathrm{d}x\,\mathrm{d}y,
\]
and so that $Q_h\bfu$ only depends on the shear-normal variable $z$. For any $\bfu\in L^2(\Omega)^3,$ the projection $Q_h\bfu$ is divergence-free and its shear-normal component vanishes everywhere, including at the free-slip boundary, but it has nonzero horizontal average unless $Q_h\bfu=\bfzero.$ Its orthogonal complement $P_h = I - Q_h$ projects onto the subspace of vector fields with streamwise and spanwise components with zero horizontal averages, but its projection does not necessarily satisfy the free-slip boundary conditions or the divergence-free condition. 

If $\bfu$ belongs to $P_{\textrm{LH}}(L^2(\Omega)^3) = H$, then it is straightforward to prove that $P_h\bfu$ still belongs to $H$ and, in fact, belongs to $\dot H.$ On the other hand, if $\bfu\in P_h(L^2(\Omega)^3),$ i.e. if $\bfu$ has zero horizontal averages of the streamwise and spanwise components, then we use the Leray-Helmholtz decomposition to write $P_{\textrm{LH}}\bfu = \bfu - \bfnabla p$ for a suitable $p$ ($p$ is periodic in $x,y$ and solves $\Delta p = \bfnabla \cdot \bfu$ with Neumann boundary condition $\partial_n p = 0$ on $z = \pm h/2$) and observe that, by periodicity, $\bfnabla p$ necessarily has zero average in the streamwise and spanwise components, so that $P_{\textrm{LH}}\bfu$ also belongs to $P_h(L^2(\Omega)^3).$ This means the two projectors $P_{\textrm{LH}}$ and $P_h$ commute with each other, and the orthogonal projector in $L^2$ onto $\dot H$ is given by $P_{\dot H} = P_{\textrm{LH}} P_h = P_h P_{\textrm{LH}}$.

We also have the decomposition $H = \dot H \oplus \dot H^{\perp_H},$ where $\dot H^{\perp_H}$ is the orthogonal complement of $\dot H$ in $H$. The space $\dot H^{\perp_H}$ consists of vector fields of the form $(\alpha(z), \beta(z), 0)$, for some square-integrable functions $\alpha, \beta$ of $z$ alone. The full decomposition of $L^2(\Omega)^3$ is then given by $L^2(\Omega)^3 = \dot H \oplus \dot H^{\perp_H} \oplus H^\perp$, and the subspace orthogonal to $\dot H$ in $L^2(\Omega)^3$ is $\dot H^\perp = \dot H^{\perp_H} \oplus H^\perp$. Notice that the body-type force $\nu U''(z)\bfe_x$ in the momentum equation \eqref{nsefluctuationprime} is of the form $(\alpha(z), 0, 0)$ and thus belongs to $\dot H^{\perp_H}$, orthogonal to $\dot H$. This term drops out when the equation is projected onto $\dot H$, as discussed in Section \ref{secmodel}.

The Stokes operator $A$ in our domain of interest is defined as $A = -P_{\dot H} \Delta$ with domain $D(A) = \{\bfu \in \dot V; A\bfu\in \dot H\} = H^2(\Omega)^3 \cap \dot V$. The Stokes operator is a positive self-adjoint operator in $\dot H$ with compact inverse.

Since $\dot H$ is a Hilbert space, we identify $\dot H$ with its dual and consider the dual $\dot V'$ of $\dot V$ with the duality product denote $\dinner{\bfg,\bfw}_{\dot V',\dot V}$, for $\bfg\in \dot V'$, $\bfw\in \dot V$, and with the operator norm 
\[ \|\bfg\|_{\dot V'} = \sup_{\|\bfnabla \bfw\|_{L^2} = 1}\dinner{\bfg,\bfw}_{\dot V',\dot V}.
\]

The Poincar\'e inequality is written as
\begin{equation}
  \label{poincareineq}
  \kappa_1^2 \|\bfu\|_{L^2}^2 \leq \|\bfnabla\bfu\|_{L^2}^2, \qquad \forall \bfu \in \dot V,
\end{equation}
where $\kappa_1$ is the smallest wavenumber of the Stokes operator, as discussed in the next two sections, and in particular in \eqref{kappazero}.

We also recall the orthogonality property
\begin{equation}
  \label{bilinortho}
  \dinner{(\bfu\cdot\bfnabla)\bfv, \bfv}_{L^2} = 0, \qquad \forall \bfu, \bfv \in \dot V,
\end{equation}
which, in particular, implies the anti-symmetry property
\begin{equation}
  \label{antisymmetry}
  \dinner{(\bfu\cdot\bfnabla)\bfv, \bfw}_{L^2} = - \dinner{(\bfu\cdot\bfnabla)\bfw, \bfv}_{L^2}, \qquad \forall \bfu, \bfv, \bfw \in \dot V,
\end{equation}

\subsection{Weak solutions and functional formulation}

Multiplying the equation \eqref{nsefluctuation} by a test function $\bfv \in \dot V$ and integrating over the domain $\Omega$ yield the weak formulation of the problem. We say that $\bfu$ is a weak solution of the system \eqref{nsefluctuation} if $\bfu \in L^\infty(0, T; \dot H) \cap L^2(0, T; \dot V)$ for all $T > 0$, and if for all $\bfv \in \dot V$ and for almost every $t > 0,$ we have
\begin{multline}
  \label{weakformulation}
  \frac{\partial}{\partial t}\dinner{\bfu, \bfv}_{L^2} + \nu \dinner{\bfnabla\bfu, \bfnabla\bfv}_{L^2} + \dinner{(\bfu\cdot\bfnabla)\bfu, \bfv}_{L^2} \\ 
  = -\dinner{U(z)\partial_x\bfu + w U'(z) \bfe_x, \bfv}_{L^2}.
\end{multline}

The functional-analytic formulation of the problem is obtained by writing the above weak formulation in terms of the Stokes operator $A$ and of the bilinear term $B(\bfu, \bfv) = P_{\dot H}((\bfu\cdot\bfnabla)\bfv)$, so that we can write
\begin{equation}
  \label{functionalformulation}
  \frac{\mathrm{d}\bfu}{\mathrm{d}t} + \nu A\bfu + B(\bfu, \bfu) = - P_{\dot H}(U(z)\partial_x\bfu) - P_{\dot H}(w U'(z) \bfe_x).
\end{equation}
Notice that $U(z)\partial_x\bfu$ and $w U'(z) \bfe_x$ have zero horizontal average but are not divergence free, so the projector $P_{\dot H}$ is needed in the right-hand side, as well.

We write \eqref{functionalformulation} more concisely as 
\begin{equation}
  \label{functionalformulationwithF}
  \frac{\mathrm{d}\bfu}{\mathrm{d}t} = \bfF(\bfu),
\end{equation}
where
\begin{equation}
  \label{functionalformulationFdef}
  \bfF(\bfu) =  - P_{\dot H} (U(z)\partial_x\bfu) - P_{\dot H}(w U'(z) \bfe_x) - \nu A\bfu - B(\bfu, \bfu).
\end{equation}
For $\bfu\in\dot V,$ we have $\bfF(\bfu)\in\dot V',$ and the weak solution satisfies $\bfu_t \in L^{4/3}(0, T, \dot V').$

Formally taking $\bfv = \bfu$ in the weak formulation \eqref{weakformulation} and using the orthogonality property \eqref{bilinortho} yield the energy balance equation
\begin{equation} 
  \label{energybalance}
  \frac{1}{2}\frac{\mathrm{d}}{\mathrm{d}t}\|\bfu\|_{L^2}^2 + \nu \|\bfnabla\bfu\|_{L^2}^2 = -\dinner{U(z)\partial_x\bfu + w U'(z) \bfe_x, \bfu}_{L^2}.
\end{equation}

\subsection{The spectral decomposition of the Stokes operator}

The Stokes operator in this periodic channel with free-slip boundary conditions and zero horizontal average of the streamwise and spanwise components possesses a complete orthonormal basis $\{\mathbf{w}_{j,l,k,\iota}\}_{j,l,k,\iota}$ of eigenfunctions in $\dot H$, with associated strictly positive eigenvalues $\lambda_{j,l,k},$ for $\iota=1,2,$ so that $A\mathbf{w}_{j,l,k,\iota}=\lambda_{j,l,k}\mathbf{w}_{j,l,k,\iota}$, for all indices in
\[
    \Ical = \{(j,l,k,\iota); \; j,l\in\ZZ, k\in\NN_0, \iota\in\{1, 2\}, j^2 + l^2 \neq 0\},
\]
where $\ZZ$ is the set of integers and $\NN_0$ is the set of nonnegative integers. The set of stricly positive integers is denoted by $\NN$.

Due to the geometry, the eigenfunctions are of the form
\begin{multline}
\label{evnew}
  \bfw_{j,l,k,\iota}(x,y,z)
    =\bfw^h_{j,l}(x,y)\odot\bfw_{j,l,k,\iota}^v(z) = \begin{pmatrix} u^h_{j,l}(x,y) \\ v^h_{j,l}(x,y) \\ w^h_{j,l}(x,y) \end{pmatrix} \odot \begin{pmatrix} u^v_{j,l,k,\iota}(z) \\ v^v_{j,l,k,\iota}(z) \\ w^v_{j,l,k,\iota}(z) \end{pmatrix} \\ 
    = \begin{pmatrix} u^h_{j,l}(x,y) u_{j,l,k,\iota}^v(z) \\ v^h_{j,l}(x,y) v_{j,l,k,\iota}^v(z) \\ w^h_{j,l}(x,y) w_{j,l,k,\iota}^v(z)\end{pmatrix},
\end{multline}
where $\odot$ denotes the Hadamard product. The components of horizontal terms are the real counterparts of the complex eigenfunctions $\exp(2\pi i(jx/L_x+ly/L_y))$, while each vertical term $\mathbf{w}_{j,l,k,\iota}^v(z)$ is given by
\begin{multline*} 
    \mathbf{w}_{j,l,k,\iota}^v(z) = \bigg(A_{j,l,k,\iota}\cos\left( \frac{k\pi}{h}\left(z + \frac{h}{2}\right)\right), B_{j,l,k,\iota}\cos\left( \frac{k\pi}{h}\left(z + \frac{h}{2}\right)\right), \\ C_{j,l,k,\iota}\sin\left( \frac{k\pi}{h}\left(z + \frac{h}{2}\right) \right) \bigg), \quad k = 0, 1, 2, \ldots,
\end{multline*}
with the coefficients $A_{j,l,k,\iota}, B_{j,l,k,\iota}, C_{j,l,k,\iota}$ being constrained by the divergence-free condition and scaled so that the eigenfunctions are normalized in $L^2(\Omega)^3$. 

Not all combinations of the real counterparts of the complex horizontal eigenfunctions are allowed in the components of $\bfw^h_{j,l}(x,y).$ Due to the zero horizontal average condition of the streamwise and spanwise components, no term with $j=l=0$ is allowed. Moreover, due to the divergence-free condition, the combinations need to match those of the vertical eigenfunctions. More precisely, the allowed combinations are the following:
\begin{align*}
    \bfw^h_{j,l}(x,y) & = \frac{1}{L_x^{1/2} L_y^{1/2}}
        \begin{pmatrix}
            \cos\left(\frac{2\pi j x}{L_x}\right)\sin\left(\frac{2\pi l y}{L_y}\right) \\
            \sin\left(\frac{2\pi j x}{L_x}\right)\cos\left(\frac{2\pi l y}{L_y}\right) \\
            \sin\left(\frac{2\pi j x}{L_x}\right)\sin\left(\frac{2\pi l y}{L_y}\right)
        \end{pmatrix}, & & j, l > 0, \\ 
    \bfw^h_{j,l}(x,y) & = \frac{1}{L_x^{1/2} L_y^{1/2}}
        \begin{pmatrix}
            \cos\left(\frac{2\pi j x}{L_x}\right)\cos\left(\frac{2\pi l y}{L_y}\right) \\
            \sin\left(\frac{2\pi j x}{L_x}\right)\sin\left(\frac{2\pi l y}{L_y}\right) \\
            \sin\left(\frac{2\pi j x}{L_x}\right)\cos\left(\frac{2\pi l y}{L_y}\right)
        \end{pmatrix}, & & j < 0, l \geq 0, \\ 
    \bfw^h_{j,l}(x,y) & = \frac{1}{L_x^{1/2} L_y^{1/2}}
        \begin{pmatrix}
            \sin\left(\frac{2\pi j x}{L_x}\right)\sin\left(\frac{2\pi l y}{L_y}\right) \\
            \cos\left(\frac{2\pi j x}{L_x}\right)\cos\left(\frac{2\pi l y}{L_y}\right) \\
            \cos\left(\frac{2\pi j x}{L_x}\right)\sin\left(\frac{2\pi l y}{L_y}\right)
        \end{pmatrix}, & & j \geq 0, l < 0, \\ 
    \bfw^h_{j,l}(x,y) & = \frac{1}{L_x^{1/2} L_y^{1/2}}
        \begin{pmatrix}
            \sin\left(\frac{2\pi j x}{L_x}\right)\cos\left(\frac{2\pi l y}{L_y}\right) \\
            \cos\left(\frac{2\pi j x}{L_x}\right)\sin\left(\frac{2\pi l y}{L_y}\right) \\
            \cos\left(\frac{2\pi j x}{L_x}\right)\cos\left(\frac{2\pi l y}{L_y}\right)
        \end{pmatrix}, & & j, l \leq 0, \;j + l < 0.
\end{align*}

The divergence-free constraint on the multiplicative factors of the vertical eigenfunctions depends on the combinations:
\begin{align*}
    & - \frac{2\pi j}{L_x} A_{j,l,k,\iota} - \frac{2\pi l}{L_y} B_{j,l,k,\iota} + \frac{k\pi}{h} C_{j,l,k,\iota} = 0, & & j, l > 0, \\ 
    & - \frac{2\pi j}{L_x} A_{j,l,k,\iota} + \frac{2\pi l}{L_y} B_{j,l,k,\iota} + \frac{k\pi}{h} C_{j,l,k,\iota} = 0, & & j < 0, l \geq 0, \\ 
    & + \frac{2\pi j}{L_x} A_{j,l,k,\iota} - \frac{2\pi l}{L_y} B_{j,l,k,\iota} + \frac{k\pi}{h} C_{j,l,k,\iota} = 0, & & j \geq 0, l < 0, \\ 
    & + \frac{2\pi j}{L_x} A_{j,l,k,\iota} + \frac{2\pi l}{L_y} B_{j,l,k,\iota} + \frac{k\pi}{h} C_{j,l,k,\iota} = 0, & & j, l \leq 0, \;j + l < 0.
\end{align*}

For each triplet of indices $(j, l, k)$, the corresponding constraint is one nontrivial equation with three unknowns $A_{j,l,k,\iota}, B_{j,l,k,\iota}, C_{j,l,k,\iota}$, so that there is a two-dimensional solution space, corresponding to $\iota = 1, 2.$ 

The eigenvalues are given by the sum of the horizontal and vertical contributions, which are decoupled due to the geometry and the boundary conditions. The horizontal contribution is given by the sum of squares of the horizontal wavenumbers, while the vertical contribution is given by the square of the vertical wavenumber. More precisely, the eigenvalues are
\begin{equation}
    \label{eigenvalues}
    \lambda_{j,l,k} = \frac{4\pi^2 j^2}{L_x^2} + \frac{4\pi^2 l^2}{L_y^2} + \frac{k^2 \pi^2}{h^2}, \quad j, l \in \ZZ, k\in \NN_0, j^2 + l^2 \neq 0.
\end{equation}
The multiplicity of each eigenvalue depends on the indices. For each triplet $(j, l, k),$ the divergence-free condition leaves two eigenfunctions. Moreover, since $\lambda_{j,l,k}$ depends only on $j^2$, $l^2$, and $k^2$, and $j$ and $l$ can be negative integers, the same eigenvalue is obtained by changing the signs of $j$ and $l$. If exactly one of $j$ or $l$ is zero, there are two sign combinations, either $\pm j$ or $\pm l,$ and the multiplicity is at least $4$. If both $j$ and $l$ are non-zero, there are four sign combinations, $(\pm j, \pm l)$, giving multiplicity at least $8$. The multiplicity can be higher if different triplets yield the same value of $\lambda_{j,l,k}$.

The smallest eigenvalue depends on the shape of the domain and can be written as
\begin{equation}
    \lambda_1 =  \frac{4\pi^2}{(\max\{L_x, L_y\})^2}.
\end{equation}
The multiplicity of the first eigenvalue is then either $4$, if $L_x \neq L_y$, or $8$, if $L_x = L_y.$

\subsection{Spectral representation of vector fields and horizontal spectral projections}
\label{secspecdecomposition}

We write the spectral expansion of a vector field $\bfu\in\dot H$ in the orthonormal basis of the eigenfunctions of the Stokes operator as
\begin{equation}
  \mathbf{u}(x,y,z) = \sum_{j,l,k,\iota} \hat u_{j,l,k, \iota} \mathbf{w}_{j,l,k,\iota}(x,y,z), 
    \qquad \hat u_{j,l,k,\iota} = \dinner{\mathbf{u},\mathbf{w}_{j,l,k,\iota}}_{L^2}.
\end{equation}
For notational simplicity, whenever we write a summation with one or more of the indices $j, l, k, \iota,$ it is assumed that each of them varies within their range in $\Ical,$ namely $j, l\in\ZZ,$ $k\in\NN_0,$ $\iota=1,2,$ up to any other restriction made explicit in the summation.

To each eigenvalue $\lambda=\lambda_{j,l,k}$ we associate a wavenumber $\kappa=\kappa_{j,l,k} = \lambda_{j,l,k}^{1/2}.$ 
The wavenumbers $\kappa_{j,l,k}$ form a discrete set in the interval $(0, \infty)$, and, depending on the value of $\kappa$, there may be no corresponding $\kappa_{j,l,k},$ and the associated component is always zero. 

Due to the nature of shear flows, we focus on the horizontal modes, separating the large and small cross-shear scales. We denote the horizontal wavenumbers by
\begin{equation} 
    \bar\kappa_{j,l} = \sqrt{\frac{4\pi^2 j^2}{L_x^2} + \frac{4\pi^2 l^2}{L_y^2}}, \qquad j, l \in \ZZ, j^2 + l^2 \neq 0.
\end{equation}
Since $k$ in \eqref{eigenvalues} is allowed to be zero, the smallest horizontal wavenumber $\bar\kappa_1$ and the smallest total wavenumber $\kappa_1$ coincide with each other and are given by
\begin{equation}
  \label{kappazero}
  \kappa_1 = \bar\kappa_1 = \sqrt{\frac{4\pi^2}{(\max\{L_x, L_y\})^2}}.
\end{equation}

The component of a vector field with a range of wavenumbers $\bar\chi$, say between $\bar\kappa_a$ and $\bar\kappa_b$, with $0 \leq \bar\kappa_a < \bar\chi \leq \bar\kappa_b \leq \infty$, is denoted by
\begin{equation}
  \bfu_{\bar\kappa_a, \bar\kappa_b}(x,y,z) = \sum_{\bar\kappa_a < \bar\chi \leq \bar\kappa_b} \sum_{\bar\kappa_{j,l} = \bar\chi}\hat u_{j,l,k,\iota} \bfw_{j,l,k,\iota}(x,y,z).
\end{equation}
The condition $\bar\kappa_{j,l} = \bar\chi$ in the summation is to be read as $\{(j,l,k,\iota)\in\Ical; \; \bar\kappa_{j,l} = \bar\chi\}.$ Moreover, while the overbar in $\bar\kappa_{j,l}$ is part of its definition and is essential to refer to the horizontal wavenumber associated with the choice of the indices, the wavenumber $\bar\chi$ is an arbitrary positive number, and its overbar is just a visual aid to remind us that we are looking at the horizontal modes. By taking $\bar\kappa_a = 0$ and $\bar\kappa_b=\infty$, we recover $\bfu_{0,\infty} = \bfu$.

A fundamental property of this decomposition is a sharper spectral version of the Poincar\'e inequality. Indeed, since $\bfu_{\bar\kappa_a,\bar\kappa_b}$ only contains modes with horizontal wavenumbers larger than $\bar\kappa_a$ and smaller than or equal to $\bar\kappa_b,$, i.e. for each eigenfunction $\bfw_{j,l,k,\iota}$ in the range, the horizontal wavenumber satisfies $\bar \kappa_a < \bar\kappa_{j,l} \leq \bar \kappa_b ,$ we have
\begin{multline}
    \label{spectralpoincare}
    \|\bfu_{\bar\kappa_a,\bar\kappa_b}\|_{L^2}^2 = \sum_{\bar\kappa_a < \bar\kappa_{j,l} \leq \bar\kappa_b} |\hat u_{j,l,k,\iota}|^2 \leq \frac{1}{\bar\kappa_a^2} \sum_{\bar\kappa_a < \bar\kappa_{j,l} \leq \bar\kappa_b} \bar\kappa_{j,l}^2 |\hat u_{j,l,k,\iota}|^2 \\ 
    = \frac{1}{\bar\kappa_a^2}\|\bfnabla_{x,y}\bfu_{\bar\kappa_a,\bar\kappa_b}\|_{L^2}^2 \leq \frac{1}{\bar\kappa_a^2}\|\bfnabla\bfu_{\bar\kappa_a,\bar\kappa_b}\|_{L^2}^2,
\end{multline}
for any $\bar\kappa_1 \leq \bar\kappa_a < \bar\kappa_b \leq \infty$, and, similarly,
\begin{equation}
    \label{reversespectralpoincare}
    \|\bfnabla_{x,y}\bfu_{\bar\kappa_a,\bar\kappa_b}\|_{L^2}^2 \leq \bar\kappa_b^2 \|\bfu_{\bar\kappa_a,\bar\kappa_b}\|_{L^2}^2,
\end{equation}
for any $0 \leq \bar\kappa_a < \bar\kappa_b < \infty.$

More significantly for the cascade is the orthogonality of the trilinear terms associated with the right hand side of \eqref{functionalformulation} when we consider different vector fields with non-overlapping ranges of horizontal wavenumbers, i.e.
\begin{equation}
    \label{horizontalorthogonality}
    \dinner{U(z)\partial_x\bfu_{\bar\kappa_a,\bar\kappa_b}, \bfv_{\bar\kappa_c,\bar\kappa_d}}_{L^2} = 0, \quad \dinner{w_{\bar\kappa_a,\bar\kappa_b} U'(z) \bfe_x, \bfv_{\bar\kappa_c,\bar\kappa_d}}_{L^2} = 0,
\end{equation}
for all $\bfu, \bfv\in V,$ and with $(\bar\kappa_a, \bar\kappa_b]\cap(\bar\kappa_c,\bar\kappa_d] = \emptyset$, where $w_{\bar\kappa_a,\bar\kappa_b}$ is the shear-normal component of $\bfu_{\bar\kappa_a,\bar\kappa_b}$. See Section \ref{secenergybudget} for more details.

\subsection{Stationary statistical solutions}
\label{secsss}

For our rigorous treatment of ensemble averages, we use the notion of statistical solution, which is a mathematical object modeling the evolution of the probability distribution of velocity fields governed by the Navier--Stokes equations \cite{foias72,vishikfursikov88,fmrt2001a,frt2010a}. Statistically stationary flows are then modeled by stationary statistical solution \cite{foias73,fmrt2001a,frt2019}.

There are essentially two notions of statistical solution in the literature, one in which the probability measures are defined directly on the phase space of velocity fields, and another in which the probability measures are defined on the space of trajectories of velocity fields and then projected onto the phase space. The second notion has better analytical properties, while the first one is more general. Since we do not require the better analytical properties of the second notion, we use the first one, which is the original notion introduced by Foias and Prodi in \cite{foias72, foias73}.

A \textbf{Foias-Prodi stationary statistical solution} of the NSE \eqref{nsefluctuation} is a Borel probability measure $\mu$ on $\dot H$ satisfying three conditions: 
\begin{enumerate}
  \item The mean enstrophy is finite, i.e.
\begin{equation}
  \label{finitemeanenstrophy}  
  \int_{\dot H} \|\bfnabla\bfu\|_{L^2}^2 \;\rmd\mu(\bfu) < \infty;
\end{equation}
  \item It satisfies the stationary statistical equation
\begin{equation} 
  \label{stationarystatisticaleq}
  \int_{\dot H} \dinner{\bfF(\bfu), \Phi'(\bfu)}_{\dot V',\dot V} \;\rmd\mu(\bfu) = 0,
\end{equation}
for any cylindrical test function $\Phi$, where $\bfF(\bfu)$ is the right hand side of the functional form of the NSE, given in \eqref{functionalformulationFdef}; $\Phi'(\bfu)$ is the Fr\'echet derivative of $\Phi$ in $H$; and a cylindrical test function is any function of the form $\Phi(\bfu) = \varphi(\dinner{\bfu,\bfv_1}_{L^2}, \ldots, \dinner{\bfu,\bfv_k}_{L^2})$, where $\bfv_j\in V$, $j=1,\ldots, k$, and $\varphi : \RR^k\rightarrow \RR$ is a $\Ccal^1$ function with compact support;
\item It satisfies the mean energy-type inequality
\begin{equation} 
  \label{stationaryenergyinequality}
  \int_{\dot H} \psi'(\|\bfu\|_{L^2}^2) \left\{ \nu\|\bfnabla\bfu\|_{L^2}^2 + \dinner{U(z)\partial_x\bfu + w U'(z) \bfe_x, \bfu}_{L^2} \right\} \;\rmd\mu(\bfu) \leq 0, 
\end{equation}
for any nonnegative, nondecreasing, continuously-differentiable real-valued function $\psi$ on $[0,\infty)$ with bounded derivative.
\end{enumerate}

The notion of stationary statistical solution extends the notion of invariant measure to three-dimensional Navier--Stokes equations, for which there is no known global well-posedness result.

\subsection{Ensemble averages and mean quantities}
\label{secensembleavg}
In the conventional theory of turbulence, one is interested in quantities such as mean kinetic energy, mean energy dissipation, mean velocity at a given point, mean velocity correlations at different points, and so on. Here we define the relevant mean quantities in the context of stationary statistical solutions, which are the rigorous mathematical objects we use to represent ensemble averages of flows in statistical equilibrium in time.

With a Borel probability measure $\mu$ as stationary statistical solution, the corresponding ensemble average of a given quantity $\varphi$ of the state of the flow is denoted by
\[ \average{\varphi} = \average{\varphi(\bfu)} = \int_{\dot H} \varphi(\bfu)\;\rmd\mu(\bfu),
\]
which is well defined for any real-valued function $\varphi \in L^1(\mu)$. In the above definition, $\bfu$ is a dummy variable both under the integration with respect to a given stationary statistical solution and in the ensemble notation $\average{\varphi} = \average{\varphi(\bfu)}$. A more common notation used in probability theory is $\mathbb{E}_{\bfu \sim \mu}[\varphi(\bfu)]$, but we prefer to use the notation $\average{\varphi(\bfu)},$ with the implicit assumption that $\bfu \sim \mu,$ as used in the conventional theory of turbulence.

Thanks to the condition \eqref{finitemeanenstrophy}, quantities such as $\varphi(\bfu) = \|\bfu\|_{L^2}^2$ and $\varphi(\bfu) = \|\bfnabla\bfu\|_{L^2}^2$ belong to $L^1(\mu)$, so that their ensemble averages are well defined. These two quantities are directly associated with the mean kinetic energy of the flow and the mean energy dissipation rate of the flow (see Section \ref{secmeanquant}).

Thanks to the H\"older and Sobolev inequalities and again to \eqref{finitemeanenstrophy}, the functionals given by 
\[ 
    \bfu \mapsto \dinner{(\bfu_{0, \bar\kappa}\cdot\bfnabla)\bfu, \bfu}_{L^2}, \quad \bfu \mapsto \dinner{(\bfu\cdot\bfnabla)\bfu, \bfu_{0, \bar\kappa}}_{L^2}, \quad \bfu \mapsto \dinner{(\bfu\cdot\bfnabla)\bfu_{0, \bar\kappa}, \bfu}_{L^2}
\]
also belong to $L^1(\mu)$, and their ensemble averages are well defined.

\subsection{Mean scalar physical quantities}
\label{secmeanquant}

Consider the ensemble average $\average{\cdot}$ with respect to a given stationary statistical solution $\mu$, as discussed in  Section \ref{secensembleavg}.

Since the fluid has constant density, we consider quantities defined per unit mass, so we do not need to carry the density of the fluid in the estimates. We take $1/\kappa_1$ as the characteristic length and, by denoting $\rho$ as the mass density of the fluid, we have that $\rho \kappa_1^{-3}$ is the characteristic unit of mass. With this characteristic scale, the \textbf{mean kinetic energy per unit mass of the flow} is given by
\[ K = \frac{1}{2}\frac{\kappa_1^3}{\rho}\int_\Omega \rho |\bfu(\bfx)|^2\;\rmd\bfx = \frac{\kappa_1^3}{2} \average{\|\bfu\|_{L^2}^2}.
\]
Similarly, the \textbf{mean kinetic energy dissipation rate per unit mass} of the flow and of its fluctuation field are given by
\begin{equation}
  \label{defepsilontildeepsilon}
  \epsilon = \nu \kappa_1^3 \average{\|\bfnabla\bfu\|_{L^2}^2}.
\end{equation}
Thanks to \eqref{finitemeanenstrophy}, these quantities are well defined. 

The strength of the shear is taken here to be
\begin{equation}
    \label{sheargradient}
    S = \|U'\|_{L^\infty}.
\end{equation}
For shear flows, a relevant quantity that enters our cascade result is the \textbf{viscous shear wavenumber}
\begin{equation}
    \label{defviscousshearkappa}
    \kappa_S = \sqrt{\frac{S}{\nu}}.
\end{equation}
The inverse $\ell_S = \kappa_S^{-1}$ is the characteristic length scale at which the corresponding time scale of viscous diffusion $\ell_S^2/\nu$ in those scales is comparable with the time scale $1/S$ of the mean shear.

Another fundamental quantity associated with shear flows is the \textbf{Corrsin wave\-number}
\begin{equation}
    \label{defcorrsinwavenumber}
    \kappa_C = \frac{S^{3/2}}{\epsilon^{1/2}}.
\end{equation}
Its inverse $\ell_C = \kappa_C^{-1} = (\epsilon/S^3)^{1/2}$ is the Corrsin length scale, introduced in \cite{corrsin1958}, which represents the scale at which the corresponding turnover time $(\ell_C^2/\epsilon)^{1/3}$ of those scale is roughly equal to the characteristic time $1/S$ of the mean shear. It marks the transition between large-scale anisotropic motion and small-scale isotropic turbulence.

The \textbf{mean kinetic energy dissipation rate of a range of horizontal modes} is denoted by
\begin{equation}
    \label{defepsilonkappaback}
    \epsilon_{\bar\kappa_a, \bar\kappa_b} = \nu \kappa_1^3 \average{\|\bfnabla \bfu_{\bar\kappa_a, \bar\kappa_b}\|_{L^2}^2},
\end{equation}
for $0 \leq \bar\kappa_a < \bar\kappa_b\leq \infty$. This is the viscous dissipation rate of the mean kinetic energy of the components of the velocity field contained in the modes between $\bar\kappa_a$ and $\bar\kappa_b$. This is one of the terms that appear in the mean energy-budget equation for those modes. The other term comes from the redistribution of kinetic energy to and from the other modes. 

The energy dissipation rate \eqref{defepsilonkappaback} can be split into horizontal and vertical parts:
\begin{equation}
    \label{defepsilonkappabacksplit}
    \epsilon_{\bar\kappa_a, \bar\kappa_b} = \epsilon_{\bar\kappa_a, \bar\kappa_b}^{x,y} + \epsilon_{\bar\kappa_a, \bar\kappa_b}^z,
\end{equation}
where
\begin{equation}
    \label{defepsilonkappabacksplittedterms}
    \epsilon_{\bar\kappa_a, \bar\kappa_b}^{x,y} = \nu \kappa_1^3 \average{\|\bfnabla_{x,y} \bfu_{\bar\kappa_a, \bar\kappa_b}\|_{L^2}^2}, \quad \epsilon_{\bar\kappa_a, \bar\kappa_b}^{z} = \nu \kappa_1^3 \average{\|\partial_z \bfu_{\bar\kappa_a, \bar\kappa_b}\|_{L^2}^2}.
\end{equation} 

For the redistribution term in the energy-budget equation, we define the \textbf{mean (low-to-high) rate of transfer of kinetic energy per unit mass of the flow}, also known as \textbf{mean (low-to-high) energy flux}, to the component of the fluctuation field with horizontal wavenumbers larger than $\bar\kappa$, for $\bar\kappa \geq \bar\kappa_1,$ by
\begin{equation}
  \label{fluxalltofluct}
  \overset{\rightarrow}{\Pi}_{\bar\kappa} = - \kappa_1^3\average{\dinner{(\bfu\cdot \bfnabla)\bfu,\bfu_{\bar\kappa,\infty}}_{L^2}} = - \kappa_1^3\average{\dinner{(\bfu\cdot \bfnabla)\bfu_{0,\bar\kappa}, \bfu_{\bar\kappa,\infty}}_{L^2}},
\end{equation}
where the last equality follows from the orthogonality property \eqref{bilinortho}.

We prove in Section \ref{secfluxasymptotic} that the mean energy flux converges, as $\bar\kappa\rightarrow\infty$, to a nonnegative \textbf{singularity flux}, or \textbf{flux at infinity},
\begin{equation}
    \label{singularityflux}
    \overset{\rightarrow}{\Pi}_{\infty} = \lim_{\bar\kappa\rightarrow\infty} \overset{\rightarrow}{\Pi}_{\bar\kappa},
\end{equation}
and, based on that, we define the \textbf{restricted (low-to-high) energy flux}
\begin{equation}
    \label{restrictedflux}
    \overset{\rightarrow}{\Pi}_{\bar\kappa}^* = \overset{\rightarrow}{\Pi}_{\bar\kappa} - \overset{\rightarrow}{\Pi}_{\infty}.
\end{equation}
The singularity flux accounts somehow to energy that is lost to any singularity of the stationary statistical solution, akin to the possible loss of regularity of individual weak solutions of the three-dimensional Navier--Stokes equations. The restricted energy flux  measures the energy that is actually transfered to the higher finite wavenumber modes, without leaking to infinity (see Section \ref{secfluxasymptotic}). If the stationary statistical solution is smooth, the flux at infinity is zero, and both fluxes at finite wavenumbers coincide. 

We also define the \textbf{mean shear energy-production rate per unit mass of a range of horizontal modes} of the flow by
\begin{equation}
    \label{defenergyproduction}
    \Pcal_{\bar\kappa_a, \bar\kappa_b} = - \kappa_1^3\average{\dinner{w U'(z) \bfe_x,\bfu_{\bar\kappa_a, \bar\kappa_b}}_{L^2}},
\end{equation}
for the modes between $\bar\kappa_a$ and $\bar\kappa_b$, for $0 \leq \bar\kappa_a < \bar\kappa_b \leq \infty,$

An energy cascade is expected to occur in the inertial range defined by $\kappa_C \ll \kappa \ll \kappa_\eta,$ where
\begin{equation}
    \kappa_\eta = \left(\frac{\epsilon}{\nu^3}\right)^{1/4}
\end{equation}
is the Kolmogorov dissipation wavenumber, associated with the Kolmogorov dissipation length scale $\ell_\eta = \kappa_\eta^{-1},$ where most of the kinetic energy dissipation takes place.

For controlling the mean energy production for the energy cascade, we finally consider the \textbf{horizontal Taylor wavenumber of a range of horizontal modes}
\begin{equation}
    \label{defKcal}
    \Kcal_{\bar\kappa_a,\bar\kappa_b} = \left(\frac{\average{\|\bfnabla_{x,y}\bfu_{\bar\kappa_a,\bar\kappa_b}\|_{L^2}^2}}{\average{\|\bfu_{\bar\kappa_a,\bar\kappa_b}\|_{L^2}^2}}\right)^{1/2}.
\end{equation}
The length scale $\Kcal_{\bar\kappa_a,\bar\kappa_b}^{-1}$ is similar in nature to the Taylor-like wavenumber defined by $\kappa_\tau = \sqrt{\average{\|\bfnabla\bfu\|_{L^2}^2}/\average{\|\bfu\|_{L^2}^2}} = \sqrt{\epsilon/2\nu K}$, except that $\Kcal_{\bar\kappa_a,\bar\kappa_b}^{-1}$ is concentrated on a given range of wavenumbers.

Thanks to the spectral estimates \eqref{spectralpoincare} and \eqref{reversespectralpoincare}, we have
\begin{equation}
    \label{Kcallowerupperbounds}
    \bar\kappa_a \leq \Kcal_{\bar\kappa_a,\bar\kappa_b} \leq \bar\kappa_b.
\end{equation}

Note that most of the quantities defined above are dependent on the given statistical solution $\mu$, as well as on the imposed shear profile $U=U(z)$, the geometric parameters $L_x$, $L_y$, and $h$, and the kinematic viscosity $\nu$ of the fluid.

By introducing the \textbf{Taylor-microscale Reynolds number} $\Rey_\tau = K/\sqrt{\nu \epsilon}$, and the \textbf{shear-rate parameter} $S^* = SK/\varepsilon$, we can relate the different characteristic wavenumbers as follows.
\[
    \frac{\kappa_C}{\kappa_S} = S \sqrt{\frac{\nu}{\varepsilon}} = \frac{S^*}{\Rey_\tau}, \qquad \frac{\kappa_S}{\kappa_\eta} = \frac{S^{1/2}\nu^{1/4}}{\epsilon^{1/4}} = \sqrt{\frac{S^*}{\Rey_\tau}}, \qquad \frac{\kappa_C}{\kappa_\eta} = \left( \frac{S^*}{\Rey_\tau} \right)^{3/2},
\]
so their relation depends on the ratio between the shear strength $S^*$ and the Taylor-microscale Reynolds number $\Rey_\tau$.

Similarly, we can also relate these quantities to the Taylor-like wavenumber $\kappa_\tau$ as follows.
\[
    \frac{\kappa_S}{\kappa_\tau} = \sqrt{2} {S^*}^{1/2}, \qquad \frac{\kappa_C}{\kappa_\tau} = \frac{\sqrt{2}{S^*}^{3/2}}{\Rey_\tau}, \qquad \frac{\kappa_\tau}{\kappa_\eta} = \frac{1}{\sqrt{2}\Rey_\tau^{1/2}}.
\]

Further relations can be obtained by introducing the \textbf{integral wavenumber} $\kappa_\ell = \epsilon/K^{3/2}$, as follows.
\[
    \frac{\kappa_C}{\kappa_\ell} = {S^*}^{3/2}, \qquad \frac{\kappa_S}{\kappa_\ell} = {S^*}^{1/2}\Rey_\tau, \qquad \frac{\kappa_\tau}{\kappa_\ell} = \frac{\Rey_\tau}{\sqrt{2}}, \qquad \frac{\kappa_\eta}{\kappa_\ell} = \Rey_\tau^{3/2}.
\]

Notice our characteristic wavenumbers are not based on isotropy and differ slightly from the classical definitions in the conventional theory of homogeneous isotropic turbulence. The difference is in the absence of numerical pre-factors, but their orders of magnitude and, especially, their relative relations, are similar.

\section{Energy-budget relations}
\label{secenergybudgets}

In this section, we let $\mu$ be a given stationary statistical solution.

\subsection{Mean energy-budget equations}
\label{secenergybudget}

By choosing 
\[ 
    \Phi(\bfu) = \frac{\kappa_1^3}{2}\psi_n(\|\bfu_{\bar\kappa_a,\bar\kappa_b}\|_{L^2}^2)
\]
in \eqref{stationarystatisticaleq}, with $0 \leq \bar\kappa_a < \bar\kappa_b <\infty$, and $\{\psi_n\}_{n\in\NN}$, with compact support and converging, uniformly on bounded sets, to the identity in $\RR$, we find, at the limit as $n\rightarrow \infty$, the mean energy budget equation with finite modes
\begin{multline}
  \label{meanenergybugdeteqraw}
  \nu \kappa_1^3\average{\|\bfnabla\bfu_{\bar\kappa_a,\bar\kappa_b}\|_{L^2}^2} + \kappa_1^3\average{\dinner{(\bfu\cdot\bfnabla)\bfu, \bfu_{\bar\kappa_a,\bar\kappa_b}}_{L^2}} \\
  = -\kappa_1^3\average{\dinner{U(z)\partial_x\bfu + w U'(z) \bfe_x,\bfu_{\bar\kappa_a,\bar\kappa_b}}_{L^2}},
\end{multline}
for $0 \leq \bar\kappa_a < \bar\kappa_b <\infty$. 

The first term on the right hand side of \eqref{meanenergybugdeteqraw} vanishes. Indeed, we first split the nonlinear term into three terms by writing $\bfu = \bfu_{0, \bar\kappa_a} + \bfu_{\bar\kappa_a,\bar\kappa_b} + \bfu_{\bar\kappa_b,\infty}$:
\begin{multline*}
    \average{\dinner{U(z)\partial_x\bfu,\bfu_{\bar\kappa_a,\bar\kappa_b}}_{L^2}} = \average{\dinner{U(z)\partial_x\bfu_{0, \bar\kappa_a},\bfu_{\bar\kappa_a,\bar\kappa_b}}_{L^2}} \\ 
     + \average{\dinner{U(z)\partial_x\bfu_{\bar\kappa_a,\bar\kappa_b},\bfu_{\bar\kappa_a,\bar\kappa_b}}_{L^2}} + \average{\dinner{U(z)\partial_x\bfu_{\bar\kappa_b, \infty},\bfu_{\bar\kappa_a,\bar\kappa_b}}_{L^2}}.
\end{multline*}
The second term on the right hand side above vanishes by integration by parts and the periodicity in the $x$ direction:
\begin{align*}
    & \average{\dinner{U(z)\partial_x\bfu_{\bar\kappa_a,\bar\kappa_b},\bfu_{\bar\kappa_a,\bar\kappa_b}}_{L^2}} \\
    & \hspace{1in} = \int_{-h/2}^{h/2} U(z) \int_0^{L_y}\int_0^{L_x} \partial_x\bfu_{\bar\kappa_a,\bar\kappa_b}(x,y,z)\cdot \bfu_{\bar\kappa_a,\bar\kappa_b}(x,y,z) \;\rmd x\,\rmd y\,\rmd z \\
    & \hspace{1in} = \int_{-h/2}^{h/2} U(z) \int_0^{L_y}\int_0^{L_x} \frac{1}{2}\partial_x|\bfu_{\bar\kappa_a,\bar\kappa_b}(x,y,z)|^2 \,\rmd x\,\rmd y\,\rmd z = 0.
\end{align*}

The other two terms vanish because the inner double integral is an $L^2$ inner product in the horizontal periodic domain between vector fields which are orthogonal in this domain because the spectral splitting is in the horizontal modes:
\begin{align*}
    & \average{\dinner{U(z)\partial_x\bfu_{0, \bar\kappa_a},\bfu_{\bar\kappa_a,\bar\kappa_b}}_{L^2}} \\
    & \qquad = \int_{-h/2}^{h/2} U(z) \int_0^{L_y}\int_0^{L_x} \partial_x\bfu_{0, \bar\kappa_a}(x,y,z)\cdot \bfu_{\bar\kappa_a,\bar\kappa_b}(x,y,z) \;\rmd x\,\rmd y\,\rmd z = 0, \\
    & \average{\dinner{U(z)\partial_x\bfu_{\bar\kappa_b, \infty},\bfu_{\bar\kappa_a,\bar\kappa_b}}_{L^2}} \\ 
    & \qquad  = \int_{-h/2}^{h/2} U(z) \int_0^{L_y}\int_0^{L_x} \partial_x\bfu_{\bar\kappa_b, \infty}(x,y,z)\cdot \bfu_{\bar\kappa_a,\bar\kappa_b}(x,y,z) \;\rmd x\,\rmd y\,\rmd z = 0.
\end{align*}

Thus, the mean energy budget equation with finite modes becomes
\begin{equation}
  \label{meanenergybugdeteq}
  \nu \kappa_1^3\average{\|\bfnabla\bfu_{\bar\kappa_a,\bar\kappa_b}\|_{L^2}^2} + \kappa_1^3\average{\dinner{(\bfu\cdot\bfnabla)\bfu, \bfu_{\bar\kappa_a,\bar\kappa_b}}_{L^2}} = -\kappa_1^3\average{\dinner{w U'(z) \bfe_x,\bfu_{\bar\kappa_a,\bar\kappa_b}}_{L^2}},
\end{equation}
for $0 \leq \bar\kappa_a < \bar\kappa_b <\infty$. 

The second term in \eqref{meanenergybugdeteq} can be written as the difference between the low-to-high rate of transfer of kinetic energy to the modes larger than or equal to $\bar\kappa_a$ and the low-to-high rate of transfer of kinetic energy to the modes larger than or equal to $\bar\kappa_b$, namely
\begin{align*}
    & \kappa_1^3\average{\dinner{(\bfu\cdot\bfnabla)\bfu, \bfu_{\bar\kappa_a,\bar\kappa_b}}_{L^2}} \\
    & \qquad = \kappa_1^3\average{\dinner{(\bfu\cdot\bfnabla)\bfu, \bfu_{\bar\kappa_a,\infty}}_{L^2}} - \kappa_1^3\average{\dinner{(\bfu\cdot\bfnabla)\bfu, \bfu_{\bar\kappa_b,\infty}}_{L^2}} \\
    & \qquad = -\overset{\rightarrow}{\Pi}_{\bar\kappa_a} + \overset{\rightarrow}{\Pi}_{\bar\kappa_b}.
\end{align*}
Using the notation \eqref{defepsilonkappaback}, \eqref{fluxalltofluct}, and \eqref{defenergyproduction}, the budget equation \eqref{meanenergybugdeteq} can be written as
\begin{equation}
  \label{meanenergybudgetineqfluctnotation}
  \epsilon_{\bar\kappa_a, \bar\kappa_b} = \overset{\rightarrow}{\Pi}_{\bar\kappa_a} - \overset{\rightarrow}{\Pi}_{\bar\kappa_b} + \Pcal_{\bar\kappa_a, \bar\kappa_b},
\end{equation}
for $0 \leq \bar\kappa_a < \bar\kappa_b <\infty$. This is a rigorous version of the statement that, in statistical equilibrium in time, the mean dissipation of energy within a certain finite range of horizontal scales balances the mean energy transfer from all the modes to that range of scales and the mean energy production within those scales.

For $\bar\kappa_b = \infty,$ we first obtain from \eqref{stationaryenergyinequality}, the mean energy inequality
\[
    \nu\kappa_1^3\average{\|\bfnabla\bfu\|_{L^2}^2} \leq - \kappa_1^3\average{\dinner{U(z)\partial_x\bfu + w U'(z) \bfe_x,\bfu}_{L^2}}.
\]
Again, the first term on the right hand side above vanishes, this time just by the same integration by parts and periodicity argument used above, namely
\begin{align*}
    \dinner{U(z)\partial_x\bfu,\bfu}_{L^2} & = \int_{-h/2}^{h/2} U(z) \int_0^{L_y}\int_0^{L_x} \partial_x\bfu(x,y,z)\cdot \bfu(x,y,z) \;\rmd x\,\rmd y\,\rmd z \\ 
    & = \frac{1}{2}\int_{-h/2}^{h/2} U(z) \int_0^{L_y}\int_0^{L_x} \partial_x|\bfu(x,y,z)|^2 \;\rmd x\,\rmd y\,\rmd z \\ 
    & = 0.
\end{align*}
Hence, we arrive at
\begin{equation}
  \label{meanenergyineq}
  \nu\kappa_1^3\average{\|\bfnabla\bfu\|_{L^2}^2} \leq - \kappa_1^3\average{\dinner{w U'(z) \bfe_x,\bfu}_{L^2}}.
\end{equation}

After subtracting \eqref{meanenergybugdeteq} with $\bar\kappa_a = 0$ and $\bar\kappa_b = \bar\kappa$ from \eqref{meanenergyineq}, we obtain
\[ 
    \nu\kappa_1^3\average{\|\bfnabla\bfu_{\bar\kappa,\infty}\|_{L^2}^2} - \kappa_1^3\average{\dinner{(\bfu\cdot\bfnabla)\bfu, \bfu_{0,\bar\kappa}}_{L^2}} \leq -\kappa_1^3\average{\dinner{w U'(z) \bfe_x,\bfu_{\bar\kappa,\infty}}_{L^2}}.
\]
Using the antisymmetry property \eqref{antisymmetry}, writing $\bfu = \bfu_{0,\bar\kappa} + \bfu_{\bar\kappa,\infty}$, and using the orthogonality property \eqref{bilinortho}, the second term becomes
\begin{align*}
    - \kappa_1^3\average{\dinner{(\bfu\cdot\bfnabla)\bfu, \bfu_{0,\bar\kappa}}_{L^2}} & = \kappa_1^3\average{\dinner{(\bfu\cdot\bfnabla)\bfu_{0,\bar\kappa}, \bfu}_{L^2}} \\ 
    & = \kappa_1^3\average{\dinner{(\bfu\cdot\bfnabla)\bfu_{0,\bar\kappa}, \bfu_{0,\bar\kappa}}_{L^2}} \\ 
    & \qquad + \kappa_1^3\average{\dinner{(\bfu\cdot\bfnabla)\bfu_{0,\bar\kappa}, \bfu_{\bar\kappa,\infty}}_{L^2}} \\
    & = \kappa_1^3\average{\dinner{(\bfu\cdot\bfnabla)\bfu_{0,\bar\kappa}, \bfu_{\bar\kappa,\infty}}_{L^2}}.
\end{align*}
Thus, we arrive at
\begin{equation}
  \label{meanenergybudgetineq}
  \nu\kappa_1^3\average{\|\bfnabla\bfu_{\bar\kappa,\infty}\|_{L^2}^2} + \kappa_1^3\average{\dinner{(\bfu\cdot\bfnabla)\bfu_{0, \bar\kappa}, \bfu_{\bar\kappa, \infty}}_{L^2}}
    \leq - \kappa_1^3\average{\dinner{w U'(z) \bfe_x,\bfu_{\bar\kappa,\infty}}_{L^2}},
\end{equation}
for $\bar\kappa \geq \bar\kappa_1$. Thanks to \eqref{horizontalorthogonality}, the right hand side of \eqref{meanenergybudgetineq} is effectively driven by $w_{\bar\kappa,\infty}$, which can still be significantly large for $\bar\kappa$ in the inertial range, but nevertheless suitably controlled by the mean energy dissipation in the large scales, as estimated in Section \ref{secestimatingproductionterm}.

Using the notation \eqref{defepsilonkappaback}, \eqref{fluxalltofluct}, and \eqref{defenergyproduction} in Section \ref{secmeanquant}, this can be written as
\begin{equation}
  \label{meanenergybudgetineqfluct}
  \epsilon_{\bar\kappa, \infty} \leq \overset{\rightarrow}{\Pi}_{\bar\kappa} + \Pcal_{\bar\kappa, \infty},
\end{equation}
for $\bar\kappa \geq \bar\kappa_1$. If it were an equality, this would be a rigorous version of the statement that, in statistical equilibrium in time, the mean dissipation of energy within a certain range of horizontal scales, balances the mean flux of energy to that range of scales and the mean energy production at those scales. The inequality comes from the possible lack of regularity of stationary statistical solutions in the three-dimensional NSE, akin to the possible lack of regularity of weak solutions of the three-dimensional NSE, and the associated possible loss of energy to modes ``at infinity'', due to the development of singularities in the flow.

\subsection{Estimating the production term}\label{secestimatingproductionterm}
We first rewrite the production term using that horizontal modes outside a certain range $[\bar\kappa_a, \bar\kappa_b)$ are orthogonal to these modes (see \eqref{horizontalorthogonality}), so that
\begin{align*}
    \Pcal_{\bar\kappa_a, \bar\kappa_b} & = - \kappa_1^3\average{\dinner{w U'(z) \bfe_x,\bfu_{\bar\kappa_a, \bar\kappa_b}}_{L^2}} \\
    & = - \kappa_1^3\average{\dinner{w U'(z) \bfe_x,u_{\bar\kappa_a, \bar\kappa_b}\bfe_x}_{L^2}} \\
    & = - \kappa_1^3\average{\dinner{w_{\bar\kappa_a, \bar\kappa_b} U'(z) \bfe_x,u_{\bar\kappa_a, \bar\kappa_b}\bfe_x}_{L^2}},
\end{align*} 
where $u_{\bar\kappa_a, \bar\kappa_b}$ and $w_{\bar\kappa_a, \bar\kappa_b}$ are the streamwise and shear-normal components of the projection $\bfu_{\bar\kappa_a, \bar\kappa_b}.$

We estimate the production term by taking the $L^\infty$ norm of the shear gradient $U'(z)$ and the $L^2$ norm in each of the two remaining terms:
\begin{equation}
    \label{productionestimate}
    \begin{aligned}
        |\Pcal_{\bar\kappa_a, \bar\kappa_b}| & \leq \kappa_1^3\|U'\|_{L^\infty}\average{\|w_{\bar\kappa_a,\bar\kappa_b}\|_{L^2}\|u_{\bar\kappa_a,\bar\kappa_b}\|_{L^2}} \\ 
        & \leq \frac{\kappa_1^3}{2}\|U'\|_{L^\infty}\average{\|w_{\bar\kappa_a,\bar\kappa_b}\|_{L^2}^2 + \|u_{\bar\kappa_a,\bar\kappa_b}\|_{L^2}^2} \\ 
        & \leq \frac{\kappa_1^3}{2}\|U'\|_{L^\infty}\average{\|\bfu_{\bar\kappa_a,\bar\kappa_b}\|_{L^2}^2}.
        \end{aligned}
\end{equation}

Since $\bfu_{\bar\kappa_a,\bar\kappa_b}$ only contains modes with horizontal wavenumbers larger than or equal to $\bar\kappa_a$, we use the inequality \eqref{spectralpoincare} to find that
\[
    \average{\|\bfu_{\bar\kappa_a,\bar\kappa_b}\|_{L^2}^2} \leq \frac{1}{\bar\kappa_a^2}\average{\|\bfnabla_{x,y}\bfu_{\bar\kappa_a,\bar\kappa_b}\|_{L^2}^2}.
\]

Combining the inequalities, we arrive at
\[
    |\Pcal_{\bar\kappa_a, \bar\kappa_b}| \leq \frac{\kappa_1^3}{2\bar\kappa_a^2}\|U'\|_{L^\infty}\average{\|\bfnabla_{x,y}\bfu_{\bar\kappa_a,\bar\kappa_b}\|_{L^2}^2}.
\]

Using \eqref{defepsilonkappabacksplittedterms} and \eqref{defviscousshearkappa}, we write this as
\begin{equation}
    \label{productionestimateinepsilon}
    |\Pcal_{\bar\kappa_a, \bar\kappa_b}| \leq \frac{\kappa_S^2}{2\bar\kappa_a^2}\epsilon_{\bar\kappa_a, \bar\kappa_b}^{x,y},
\end{equation}
for $\bar\kappa_1 \leq \bar\kappa_a < \bar\kappa_b \leq \infty$. 

This shows that the mean energy production term within a certain range of scales is controlled by the mean energy dissipation on the same range. Since $\epsilon_{\bar\kappa_a, \bar\kappa_b}$ is bounded by the full energy dissipation $\epsilon$, we see that the mean energy production beyond a certain sufficiently large horizontal wavenumber $\bar\kappa_a$ becomes negligible compared to the mean energy dissipation.

The estimate above is good for balancing with the energy flux and the energy dissipation in \eqref{meanenergybudgetineqfluctnotation} but is not as sharp for the cascade. Instead, we use the wavenumber $\Kcal_{\bar\kappa_a,\bar\kappa_b}$ defined in \eqref{defKcal} and write
\begin{multline*}
    \average{\|\bfu_{\bar\kappa_a,\bar\kappa_b}\|_{L^2}^2} = \frac{\average{\|\bfu_{\bar\kappa_a,\bar\kappa_b}\|_{L^2}^2}}{\average{\|\bfnabla_{x,y}\bfu_{\bar\kappa_a,\bar\kappa_b}\|_{L^2}^2}}\average{\|\bfnabla_{x,y}\bfu_{\bar\kappa_a,\bar\kappa_b}\|_{L^2}^2} \\ 
    = \frac{1}{\Kcal_{\bar\kappa_a,\bar\kappa_b}^2}\average{\|\bfnabla_{x,y}\bfu_{\bar\kappa_a,\bar\kappa_b}\|_{L^2}^2} = \frac{1}{\nu\kappa_1^3\Kcal_{\bar\kappa_a,\bar\kappa_b}^2}\epsilon_{\bar\kappa_a, \bar\kappa_b}^{x,y}.
\end{multline*}
Then, using also \eqref{defviscousshearkappa}, estimate \eqref{productionestimate} yields
\begin{equation}
    \label{productionestimatewithF}
    |\Pcal_{\bar\kappa_a, \bar\kappa_b}| \leq \frac{\kappa_S^2}{2\Kcal_{\bar\kappa_a,\bar\kappa_b}^2}\epsilon_{\bar\kappa_a, \bar\kappa_b}^{x,y},
\end{equation}
for $\bar\kappa_1 \leq \bar\kappa_a < \bar\kappa_b \leq \infty$. 

\subsection{Asymptotic behavior of the flux term as the wavenumber increases}
\label{secfluxasymptotic}

From the energy-budget equation \eqref{meanenergybudgetineqfluctnotation} and using the estimate \eqref{productionestimateinepsilon} for the production term we find that
\[
    \overset{\rightarrow}{\Pi}_{\bar\kappa_a} - \overset{\rightarrow}{\Pi}_{\bar\kappa_b} = \epsilon_{\bar\kappa_a, \bar\kappa_b} - \Pcal_{\bar\kappa_a, \bar\kappa_b} \geq \epsilon_{\bar\kappa_a, \bar\kappa_b} - \frac{\kappa_S^2}{2\bar\kappa_a^2}\epsilon_{\bar\kappa_a, \bar\kappa_b}^{x,y} = \left( 1 - \frac{\kappa_S^2}{2\bar\kappa_a^2} \right)\epsilon_{\bar\kappa_a, \bar\kappa_b}^{x,y},
\]
for $\bar\kappa_1 \leq \bar\kappa_a < \bar\kappa_b <\infty$. This implies that
\begin{equation}
    \label{fluxmonotonic}
    \overset{\rightarrow}{\Pi}_{\bar\kappa_a} \geq \overset{\rightarrow}{\Pi}_{\bar\kappa_b}, \quad \forall \bar\kappa_b^2 > \bar\kappa_a^2 > \frac{\kappa_S^2}{2}.
\end{equation}
Moreover, from \eqref{meanenergybudgetineqfluct}, we find that
\begin{equation}
    \label{fluxnonnegative}
    \overset{\rightarrow}{\Pi}_{\bar\kappa} \geq \epsilon_{\bar\kappa, \infty} - \Pcal_{\bar\kappa, \infty} \geq \epsilon_{\bar\kappa, \infty} - \frac{\kappa_S^2}{2\bar\kappa^2}\epsilon_{\bar\kappa, \infty}^{x,y} \geq \left( 1 - \frac{\kappa_S^2}{2\bar\kappa^2} \right) \epsilon_{\bar\kappa, \infty}^{x,y}.
\end{equation}

Estimates \eqref{fluxmonotonic} and \eqref{fluxnonnegative} mean that, for sufficiently large wavenumber $\bar\kappa,$ namely satisfying $\bar\kappa^2 > \kappa_S^2/2$, the flux $\overset{\rightarrow}{\Pi}_{\bar\kappa}$ is nonnegative and monotonically non-increasing. In particular, this means that the asymptotic limit exists and is nonnegative:
\begin{equation}
    \label{fluxatinfinity}
    \overset{\rightarrow}{\Pi}_{\infty} = \lim_{\bar\kappa \rightarrow \infty}\overset{\rightarrow}{\Pi}_{\bar\kappa} \geq 0.
\end{equation}
This gives the ``singularity'' flux defined in \eqref{singularityflux}.

For a regular stationary statistical solution (associated e.g. with a global regular solution of the three-dimensional Navier--Stokes equations), the equality in \eqref{meanenergybudgetineqfluctnotation} holds up to $\bar\kappa_b = \infty$, and the limit $\overset{\rightarrow}{\Pi}_{\infty}$ is zero. More precisely, it suffices that a stationary statistical solution $\mu$ be carried by a bounded set in $L^4(\Omega)^3$. This, together with the finite mean enstrophy condition \eqref{finitemeanenstrophy}, allows us to apply the Lebesgues Dominated Convergence Theorem and deduce, using the ortogonality property \eqref{bilinortho}, that 
\begin{multline*} 
    \overset{\rightarrow}{\Pi}_{\infty} = \lim_{\bar\kappa\rightarrow \infty} \overset{\rightarrow}{\Pi}_{\bar\kappa} = - \kappa_1^3\lim_{\bar\kappa\rightarrow \infty}\average{\dinner{(\bfu\cdot \bfnabla)\bfu,\bfu_{\bar\kappa,\infty}}_{L^2}} \\ 
    = \kappa_1^3\lim_{\bar\kappa\rightarrow \infty}\average{\dinner{(\bfu\cdot \bfnabla)\bfu,\bfu_{0, \bar\kappa}}_{L^2}} = \kappa_1^3\average{\dinner{(\bfu\cdot \bfnabla)\bfu,\bfu}_{L^2}} = 0.
\end{multline*}
If the stationary statistical solution is not sufficiently regular, we might have a strictly positive limit $\overset{\rightarrow}{\Pi}_{\infty} > 0$, meaning the flux of energy ``leaks'' to infinity, i.e. to singularity in the stationary statistical solutions (See \cite{fmrt2001c} for a similar discussion in the fully periodic case, where this idea was first introduced).

The singularity flux is an ensemble spectral analogue of the local-energy-balance defect measure introduced by Caffarelli--Kohn--Nirenberg (CKN) in the theory of suitable weak solutions \cite{CKN1982}. Whereas the CKN defect measure captures anomalous energy dissipation locally in physical space for a given flow, the singularity flux quantifies the anomalous loss of energy escaping to infinity in spectral space and with respect to an ensemble of flows. Nevertheless, the singularity flux vanishes if and only if the strict energy equality holds almost surely (with respect to the stationary statistical solution), which in turn is equivalent to the defect measure vanishing almost surely as well.

\subsection{A restricted energy-budget equation}

Using that the energy flux converges as $\bar\kappa\rightarrow\infty,$ as proved in Section \ref{secfluxasymptotic}, and using definition \eqref{singularityflux}, we take the limit $\bar\kappa_b \rightarrow \infty$ in \eqref{meanenergybudgetineqfluctnotation} to find that
\begin{equation}
  \label{meanenergybudgetineqfluctlimit}
  \epsilon_{\bar\kappa, \infty} = \overset{\rightarrow}{\Pi}_{\bar\kappa} - \overset{\rightarrow}{\Pi}_{\infty} + \Pcal_{\bar\kappa, \infty},
\end{equation}
for $\bar\kappa_1\leq \bar\kappa <\infty$.

Using that the limit is nonnegative, we recover \eqref{meanenergybudgetineqfluct}. On the other hand, using the restricted energy flux \eqref{restrictedflux}, we write \eqref{meanenergybudgetineqfluctlimit} as the \emph{restricted energy budget equation}
\begin{equation}
  \label{meanenergybudgeteqrestrited}
  \epsilon_{\bar\kappa, \infty} = \overset{\rightarrow}{\Pi}_{\bar\kappa}^* + \Pcal_{\bar\kappa, \infty},
\end{equation}
for $\bar\kappa_1\leq \bar\kappa <\infty$.

\section{Energy cascade}
\label{secenergycascade}

Now we address our main results. From \eqref{meanenergybudgetineqfluct} and using \eqref{defepsilonkappabacksplit} and \eqref{productionestimatewithF} we find that
\begin{multline}
    \label{energyfluxlowerboundestimate1}
    \overset{\rightarrow}{\Pi}_{\bar\kappa} \geq \epsilon_{\bar\kappa, \infty} - \Pcal_{\bar\kappa, \infty} \geq \epsilon_{\bar\kappa, \infty} - \frac{\kappa_S^2}{2\Kcal_{\bar\kappa,\infty}^2}\epsilon_{\bar\kappa, \infty}^{x,y} = \epsilon_{\bar\kappa, \infty}^z + \epsilon_{\bar\kappa, \infty}^{x,y} - \frac{\kappa_S^2}{2\Kcal_{\bar\kappa,\infty}^2}\epsilon_{\bar\kappa, \infty}^{x,y}
    \\ = \epsilon_{\bar\kappa, \infty}^z + \left(1 - \frac{\kappa_S^2}{2\Kcal_{\bar\kappa,\infty}^2}\right)\epsilon_{\bar\kappa, \infty}^{x,y},
\end{multline}
for all $\bar\kappa \geq \bar\kappa_1$. The bound in \eqref{energyfluxlowerboundestimate1} can be further estimated as
\begin{multline}
    \label{energyfluxlowerboundestimate2}
    \overset{\rightarrow}{\Pi}_{\bar\kappa} \geq \left(1 - \frac{\kappa_S^2}{2\Kcal_{\bar\kappa,\infty}^2}\right)\epsilon_{\bar\kappa, \infty}^z + \left(1 - \frac{\kappa_S^2}{2\Kcal_{\bar\kappa,\infty}^2}\right)\epsilon_{\bar\kappa, \infty}^{x,y} = \left(1 - \frac{\kappa_S^2}{2\Kcal_{\bar\kappa,\infty}^2}\right)\frac{\epsilon_{\bar\kappa, \infty}}{\epsilon}\epsilon \\
    = \left(1 - \frac{\kappa_S^2}{2\Kcal_{\bar\kappa,\infty}^2}\right)\frac{\epsilon - \epsilon_{0, \bar\kappa}}{\epsilon}\epsilon = \left(1 - \frac{\kappa_S^2}{2\Kcal_{\bar\kappa,\infty}^2}\right)\left(1 - \frac{\epsilon_{0, \bar\kappa}}{\epsilon}\right)\epsilon,
\end{multline}
for all $\bar\kappa \geq \bar\kappa_1$. Then, for a range of wavenumber $\bar\kappa,$ if the range exists, we have a lower bound on the flux that is of the order of the energy dissipation rate. We state this more precisely as follows.
\begin{thm}
    \label{thmenergyfluxlowerbound}
    Let $0 < \delta \leq 1$. Then, for any $\bar\kappa$ such that 
    \begin{equation}
        \frac{\kappa_S^2}{2\Kcal_{\bar\kappa,\infty}^2} \leq \delta \quad \text{and} \quad \frac{\epsilon_{0, \bar\kappa}}{\epsilon} \leq \delta,
    \end{equation}    
    we have the following lower bound on the energy flux:
    \begin{equation}
        \overset{\rightarrow}{\Pi}_{\bar\kappa} \geq \left(1 - \delta\right)^2\epsilon.
    \end{equation}
\end{thm}

\begin{rmk}
    \label{rmklowerboundcascade}
    Theorem \ref{thmenergyfluxlowerbound} can be stated informally in the following way. For $\bar\kappa$ sufficiently large but not so large such that
    \begin{equation}
        \label{conditionKcalkandepsilonk}
        \frac{\kappa_S^2}{2} \ll \Kcal_{\bar\kappa,\infty}^2 \quad \text{and} \quad \epsilon_{0, \bar\kappa} \ll \epsilon,
    \end{equation}
    we find that
    \begin{equation}
        \label{cascadelowerboundsim}
        \overset{\rightarrow}{\Pi}_{\bar\kappa} \gtrapprox \epsilon.
    \end{equation}
    This means that within a range of wavenumbers such that the conditions \eqref{conditionKcalkandepsilonk} on $\bar\kappa$ hold, we have an energy cascade, where the flux is nearly equal to the energy dissipation rate
\end{rmk}

The same lower bound estimate obtained in \eqref{energyfluxlowerboundestimate1} for the energy flux also holds for the restricted energy flux \eqref{restrictedflux}. Indeed, using \eqref{meanenergybudgeteqrestrited}, we see that the only difference from \eqref{energyfluxlowerboundestimate1} is that the first step is an equality, 
\begin{multline}
    \label{restrictedenergyfluxlowerboundestimate1}
    \overset{\rightarrow}{\Pi}_{\bar\kappa}^* = \epsilon_{\bar\kappa, \infty} - \Pcal_{\bar\kappa, \infty} \geq \epsilon_{\bar\kappa, \infty}^z + \left(1 - \frac{\kappa_S^2}{2\Kcal_{\bar\kappa,\infty}^2}\right)\epsilon_{\bar\kappa, \infty}^{x,y} \geq \left(1 - \frac{\kappa_S^2}{2\Kcal_{\bar\kappa,\infty}^2}\right)\left(1 - \frac{\epsilon_{0, \bar\kappa}}{\epsilon}\right)\epsilon,
\end{multline}
$\bar\kappa \geq \bar\kappa_1$. Since it is an equality, we also obtain an upper bound:
\begin{multline}
    \label{restrictedenergyfluxupperboundestimate}
    \overset{\rightarrow}{\Pi}_{\bar\kappa}^* = \epsilon_{\bar\kappa, \infty} - \Pcal_{\bar\kappa, \infty} \leq \epsilon_{\bar\kappa, \infty} + \frac{\kappa_S^2}{2\Kcal_{\bar\kappa,\infty}^2}\epsilon_{\bar\kappa, \infty}^{x,y}
    \\ \leq \left(1 + \frac{\kappa_S^2}{2\Kcal_{\bar\kappa,\infty}^2}\right)\epsilon_{\bar\kappa, \infty} \leq \left(1 + \frac{\kappa_S^2}{2\Kcal_{\bar\kappa,\infty}^2}\right)\epsilon.
\end{multline}
Thus, we obtain the following result.
\begin{thm}
    \label{thmrestrictedenergyfluxbounds}
    Let $\delta > 0$. Then, for any $\bar\kappa$ such that
    \begin{equation}
        \frac{\kappa_S^2}{2\Kcal_{\bar\kappa,\infty}^2} \leq \delta \quad \text{and} \quad \frac{\epsilon_{0, \bar\kappa}}{\epsilon} \leq \delta,
    \end{equation}     
    we have the following bounds on the restricted energy flux:
    \begin{equation}
        \left(1 - \delta\right)^2\epsilon \leq \overset{\rightarrow}{\Pi}_{\bar\kappa}^* \leq \left(1 + \delta\right)\epsilon.
    \end{equation}
\end{thm}

Notice that, since the flux at infinity is nonnegative, as proved in \eqref{fluxatinfinity}, and hence the restricted energy flux is bounded by the energy flux, the result in Theorem \ref{thmenergyfluxlowerbound} actually follows from Theorem \ref{thmrestrictedenergyfluxbounds}, but we prefer to deduce the former separately to make it clear that that result does not depend on the use of the flux at infinity or the restricted energy flux.

\begin{rmk}
    Following the Remark \ref{rmklowerboundcascade}, we obtain, within the same range \eqref{conditionKcalkandepsilonk}, that the following cascade occurs for the restricted energy flux: $\overset{\rightarrow}{\Pi}_{\bar\kappa}^* \approx \epsilon.$
\end{rmk}

\section{Concluding remarks}
\label{secconcludingrmks}

Here we present some final considerations concerning our work.

\subsection{On the inertial range}
\label{concludingrmkinertialrange}

We have introduced the wavenumbers $\Kcal_{\bar\kappa,\infty}^2$ in \eqref{defKcal} and have proved that, for $\bar\kappa$ within a range given by \eqref{conditionKcalkandepsilonk}, the energy flux and the restricted energy flux satisfy the energy cascade relations
\[ 
    \overset{\rightarrow}{\Pi}_{\bar\kappa} \gtrapprox \epsilon \quad \text{and} \quad \overset{\rightarrow}{\Pi}_{\bar\kappa}^* \approx \epsilon.
\]
It is not clear, though, what exactly is the range $\bar\kappa$ in which the conditions \eqref{conditionKcalkandepsilonk} are satisfied. If we assume the Kolmogorov energy spectrum $C_K \epsilon^{2/3}\kappa^{-5/3}$ in an inertial range below the Kolmogorov wavenumber $\kappa_\eta,$ then
\[
    \average{\|\bfnabla_{x,y}\bfu_{\bar\kappa,\infty}\|_{L^2}^2} \sim C_K\epsilon^{2/3}\int_{\bar\kappa}^{\kappa_\eta} \kappa^{1/3}\;\rmd\kappa \approx C_K\epsilon^{2/3}\kappa_\eta^{4/3},
\]
and
\[
    \average{\|\bfu_{\bar\kappa,\infty}\|_{L^2}^2} \sim C_K\epsilon^{2/3}\int_{\bar\kappa}^{\kappa_\eta} \kappa^{-5/3}\;\rmd\kappa \approx C_K\epsilon^{2/3} \frac{1}{\bar\kappa^{2/3}}.
\]
Therefore,
\[ 
    \Kcal_{\bar\kappa,\infty}^2 \sim \kappa_\eta^{4/3}\bar\kappa^{2/3},
\]
and the first condition in \eqref{conditionKcalkandepsilonk} becomes
\[
    \bar\kappa \gg \frac{\kappa_S^3}{\kappa_\eta^2}.
\]
This ratio is precisely the Corrsin scale:
\[
    \frac{\kappa_S^3}{\kappa_\eta^2} = \frac{(S/\nu)^{3/2}}{(\epsilon/\nu^3)^{1/2}} = \frac{S^{3/2}}{\epsilon^{1/2}} = \kappa_C.
\]
This shows that, in contrast to the cases of body-forced turbulence, where the lower edge of the inertial range is dictated by the forcing scale, our implicit lower bound is genuinely driven by the shear scale.

Moreover, one expects most of the dissipation to occur near the Kolmogorov dissipation scales, so that
\[ \epsilon_{0, \bar\kappa} \ll \epsilon,
\]
for $\bar\kappa \ll \kappa_\eta.$ Thus, the second condition in \eqref{conditionKcalkandepsilonk} is expected to be satisfied for
\[
    \bar\kappa \ll \kappa_\eta.
\]
Therefore, assuming the Kolmogorov spectrum, the rigorous conditions \eqref{conditionKcalkandepsilonk} turn out to be precisely the heuristic conditions
\[
    \kappa_C \ll \bar\kappa \ll \kappa_\eta.
\]
This means that the conditions \eqref{conditionKcalkandepsilonk} seem to be sharp rigorous conditions for the energy cascade. With that in mind, it should be interesting to study and measure how exactly $\Kcal_{\bar\kappa,\infty}$ varies with $\bar\kappa,$ in numerical simulations of turbulent shear flows; using a pseudo-spectral method with a Fourier expansion in the horizontal cross-shear directions, which is efficient in our free-slip geometry for the fluctuation field, the computation of $\Kcal_{\bar\kappa,\infty}$ is straightforward.

For comparison, Table \ref{tablengthscales} displays the values obtained for the relevant scales in two laboratory experiments \cite{Champagne1970,Tavoularis1981} and two DNS simulations \cite{RogersMoinReynolds1986,SekimotoDongJimenez2016} of different shear flows. Champagne et al. \cite{Champagne1970} performed wind-tunnel experiments with a test section roughly one foot high and with a shear generator consisting of a grid of parallel rods with non-uniform spacing and varying resistance to create a linear mean velocity gradient ($S = 12.9\texttt{s}^{-1}$, corresponding to $S^* = 2.1$). Tavoularis and Corrsin \cite{Tavoularis1981} designed a similar experiment, except with a higher test section and an improved shear generator taylored for a more precise and higher-shear linear profile ($S = 46.8\texttt{s}^{-1}$, corresponding to $S^* = 10.5$). Rogers et al. \cite{RogersMoinReynolds1986} provided the first high-resolution DNS of a turbulent shear flow, utilizing a $128^3$ numerical grid and implementing deforming periodic boundary conditions to maintain a perfectly uniform mean shear gradient in a periodic box ($S^* \approx 5.5$). Sekimoto et al. \cite{SekimotoDongJimenez2016} conducted several DNS simulations with varying resolutions going up to $2048\times 2048\times 1024$, also maintaining a perfectly uniform shear rate ($S^* \approx 6.0$) and with a much higher Taylor-microscale Reynolds number.

\begin{table}[h]
    \centering
    \begin{tabular}{lcccc}
        \hline\hline
        Scale & \cite{Champagne1970} & \cite{Tavoularis1981} & \cite{RogersMoinReynolds1986} & \cite{SekimotoDongJimenez2016} \\ \hline
        $\ell_\eta$ & $0.22\,\texttt{mm}$ & $0.177\,\texttt{mm}$ & $\approx 0.056\ell_T$ & $\approx 0.032\ell_T$ \\
        $\ell_S$ & $1.08\,\texttt{mm}$ & $0.57\,\texttt{mm}$  & $\approx 0.135\ell_T$ & $\approx 0.129\ell_T$ \\
        $\ell_C$ & $25.2\,\texttt{mm}$ & $5.78\,\texttt{mm}$  & $\approx 0.78\ell_T$ & $\approx 2.08\ell_T$ \\
        $\ell_T$ & $5.0\,\texttt{mm}$ & $5.8\,\texttt{mm}$ & $\ell_T$ & $\ell_T$ \\
        $L$ & $\approx 50\,\texttt{mm}$ & $57\,\texttt{mm}$ & $\approx 5.5\ell_T$ & $\approx 16.7\ell_T$ \\ 
        $L/\ell_\eta$ & $\approx 230$ & $\approx 320$ & $\approx 98$ & $\approx 520$  \\ \hline\hline
    \end{tabular}
    \caption{Comparison of characteristic length scales in homogeneous shear flow: Kolmogorov scale ($\ell_\eta$), viscous shear scale ($\ell_S$), Corrsin scale ($\ell_C$), Taylor microscale ($\ell_T$), and integral scale ($L$). The values were calculated from the data reported in each of the works, with the experimental results in \cite{Champagne1970,Tavoularis1981} reported in milimeters (\texttt{mm}) and the DNS simulations in \cite{RogersMoinReynolds1986,SekimotoDongJimenez2016} being nondimensional and stated relative to the Taylor microscale $\ell_T$ in their setup. The values mentioned in this table were calculated based on their definitions, which differ slightly from our definitions by a order one pre-factor, as mentioned at the end of Section \ref{secmeanquant}.}
    \label{tablengthscales}
\end{table}

Notice that the Taylor microscale $\ell_T$ is comparable with the Corrsin scale in some of the flows, so that bounding the inertial-range scales by $\ell_T$, as done in the periodic case, is not physically realistic for shear flows, supporting the need to focus on $\Kcal_{\bar\kappa,\infty}$ and $\kappa_S$.

If the energy spectrum within the inertial range is not exactly $\kappa^{-5/3}$ but instead is proportional to $\kappa^{-\alpha}$, for some $1 < \alpha < 3$, e.g. as in the She--Leveque theory \cite{SheLeveque1994}, then 
\[ 
    \average{\|\bfnabla_{x,y}\bfu_{\bar\kappa,\infty}\|_{L^2}^2} \propto \kappa_\eta^{3 - \alpha}, \qquad \average{\|\bfu_{\bar\kappa,\infty}\|_{L^2}^2} \propto \bar\kappa^{1 - \alpha},
\] 
so that 
\[ 
    \Kcal_{\bar\kappa,\infty}^2 \sim \kappa_\eta^{3 - \alpha}\bar\kappa^{\alpha - 1},
\]
and the first condition in \eqref{conditionKcalkandepsilonk} leads to
\[
    \bar\kappa \gg \frac{\kappa_S^{2/(\alpha - 1)}}{\kappa_\eta^{(3 - \alpha)/(\alpha - 1)}} = \frac{(S/\nu)^{1/(\alpha - 1)}}{(\epsilon/\nu^3)^{(3 - \alpha)/(4\alpha - 4)}} = \kappa_C\left(\frac{\kappa_\eta}{\kappa_C}\right)^{(3\alpha - 5)/(3\alpha - 3)},
\]
showing how the lower bound varies with the deviation of $\alpha$ from $5/3$.

\subsection{On the horizontal wavenumber splitting}

We have split the vector fields by separating only their horizontal wavenumbers. The reasons are the properties \eqref{horizontalorthogonality} that have two consequences: (i) together with the periodicity condition, they make the first term in the budget equation \eqref{meanenergybugdeteqraw} vanish, and (ii) they make the part of the second term in \eqref{meanenergybugdeteqraw} with the components of $w$ with modes outside the range $[\bar\kappa_a,\bar\kappa_b)$ also vanish. This helps considerably the analysis of the energy budget equation, giving us suitable estimates to obtain conditions for the energy cascade. In one sense, this yields a weaker result since the energy flux for which we obtain the cascade only takes into account the flux of energy among horizontal scales. On the other hand, this is a stronger result since it captures the horizontal-wavenumber cascade in a possibly larger range, where the flow might not be fully isotropic.

Another analysis can be made by splitting the velocity field with respect to the full wavenumber decomposition, aiming to capture the cascade in a fully isotropic regime. This requires different estimates for the production term and also requires estimating the advection term, which does not vanish from the energy-budget equation in this case. The consequence is that we need four conditions to define the range for the full-scale energy cascade, leading to a possibly smaller inertial range, depending on the shear gradient. This result complements the current one and will be presented elsewhere. 

\subsection{Comparison with previous results based on the Taylor wavenumber}
\label{secconcludingcomparison}

The previous energy cascade results in the case of periodic boundary conditions with a body-type forcing yield the energy cascade under the condition that $\kappa_f \leq \kappa \ll \kappa_\tau,$ where $\kappa_f$ is (essentially) the largest excited mode in the forcing term and $\kappa_\tau = \average{\|\bfnabla \bfu\|_{L^2}^2}^{1/2} / \average{\|\bfu\|_{L^2}^2}^{1/2}$ is a Taylor-like wavenumber. This is based on bounding the corresponding wavenumbers $\average{\|\bfnabla\bfu_{\kappa,\infty}\|_{L^2}^2}^{1/2} / \average{\|\bfu_{\kappa,\infty}\|_{L^2}^2}^{1/2}$ (analogous to our $\Kcal_{\bar\kappa,\infty}$) from below by $\kappa$ (as in \eqref{Kcallowerupperbounds}).

We could have done the same here, but in shear flows the energy cascade might occur above $\kappa_\tau$, if the shear is strong (see Table \ref{tablengthscales}), so the estimate would not be realistic. This is the reason for working with $\Kcal_{\bar\kappa,\infty}$, here. It turns out this is a major improvement, since it seems to correspond to the expected inertial range in the conventional theory, as discussed in Section \ref{concludingrmkinertialrange}. The same improvement can be done for the periodic case, at the cost of not having an explicit condition for the inertial range.

\subsection{On the horizontal gradient}
In the definition \eqref{defKcal} of $\Kcal_{\bar\kappa_a,\bar\kappa_b},$  we have used the horizontal gradient $\bfnabla_{x,y}$. This is not strictly necessary. A similar result holds for $\Kcal_{\bar\kappa_a,\bar\kappa_b}$ defined as $\average{\|\bfnabla\bfu_{\bar\kappa_a,\bar\kappa_b}\|_{L^2}^2}^{1/2} /\average{\|\bfu_{\bar\kappa_a,\bar\kappa_b}\|_{L^2}^2}^{1/2},$ with the full gradient. We used the horizontal gradient to be consistent with the choice of splitting with horizontal wavenumbers. If we use the full gradient, we lose the bounds in \eqref{Kcallowerupperbounds}, but those are not strictly necessary.

\end{document}